\documentclass[a4paper]{article}

\usepackage{amsmath,amssymb,amsthm}
\usepackage{graphicx}
\usepackage{pgf,tikz}
\usetikzlibrary{calc,arrows,shapes,shadows,decorations.pathreplacing}
\tikzset{tip/.style = {-latex}}
\usepackage{subfig, caption}

\usepackage{hyperref}

\theoremstyle{plain}
\newtheorem{theorem}{Theorem}

\def\RR{\mathbb{R}}
\def\ZZ{\mathbb{Z}}
\def\NN{\mathbb{N}}
\def\bP{\mathbf{P}}
\def\bp{\mathbf{p}}
\def\br{\mathbf{r}}

\def\ci{\mathrm{i}}

\def\dx{\mathrm{d}x}
\def\dy{\mathrm{d}y}
\def\dt{\mathrm{d}t}

\def\ssfourpoints(#1,#2){
\pgfmathtruncatemacro{\lasti}{#1 + 1}
\pgfmathtruncatemacro{\lastii}{#1 + 2}
\coordinate (p\lasti) at (p1);
\coordinate (p\lastii) at (p2);
\coordinate (p0) at (p#1);

\foreach \h in {1,...,#1} {%
  \pgfmathtruncatemacro{\p}{\h-1}%
  \pgfmathtruncatemacro{\j}{\h+1}%
  \pgfmathtruncatemacro{\k}{\h+2}%

  \pgfmathtruncatemacro{\odd}{2*\h-1}%
  \pgfmathtruncatemacro{\even}{2*\h}%
  \coordinate (n\odd) at (p\h);
  \coordinate (n\even) at (barycentric cs:p\p=-1/16,p\h=9/16,p\j=9/16,p\k=-1/16);
}

\pgfmathsetmacro{\i}{2*#1}
\foreach \j in {1,...,\i} \coordinate (p\j) at (n\j);

\pgfmathsetmacro{\k}{#2 - 1}
\ifnum0<\k
\ssfourpoints(\i,\k)
\fi
}

\def\ssBS(#1,#2){
\pgfmathtruncatemacro{\lasti}{#1 + 1}
\coordinate (p\lasti) at (p1);
\coordinate (p0) at (p#1);

\foreach \h in {1,...,#1} {%
  \pgfmathtruncatemacro{\p}{\h-1}%
  \pgfmathtruncatemacro{\j}{\h+1}%

  \pgfmathtruncatemacro{\odd}{2*\h-1}%
  \pgfmathtruncatemacro{\even}{2*\h}%
  \coordinate (n\odd) at (barycentric cs:p\p=1/8,p\h=3/4,p\j=1/8);
  \coordinate (n\even) at (barycentric cs:p\h=1/2,p\j=1/2);
}

\pgfmathsetmacro{\i}{2*#1}
\foreach \j in {1,...,\i} \coordinate (p\j) at (n\j);

\pgfmathsetmacro{\k}{#2 - 1}
\ifnum0<\k
\ssBS(\i,\k)
\fi
}

\def\localssBS{
  \coordinate (n1) at (barycentric cs:p1=1/2,p2=1/2);
  \coordinate (n2) at (barycentric cs:p1=1/8,p2=3/4,p3=1/8);
  \coordinate (n3) at (barycentric cs:p2=1/2,p3=1/2);

\foreach \j in {1,2,3} \coordinate (p\j) at (n\j);
}

\newtheorem{propos}[theorem]{Proposition}
\newtheorem{remark}{Remark}[section]

\author{\textbf{First Author}  \and \textbf{Second Author}  \and \textbf{Third Author} \and \textbf{Fourth Author} }

\title{ \textbf{Subdivision based snakes for contour detection}}
\author{Rafael D\'{i}az Fuentes  \and Javier Pino Torres \and Victoria Hern\'andez Mederos \and Jorge Estrada Sarlabous %
\\[1em]
\{\text{rafaeldf, javier, vicky, jestrada}\}@icimaf.cu}

\date{Instituto de Cibern\'etica, Matem\'atica y F\'{i}sica }

\begin{document}
\maketitle


\begin{abstract}
In this paper we propose a method for computing the contour of an object in an image using a snake represented as a subdivision curve. The evolution of the snake is driven by its control points which are computed minimizing an energy that pushes the snake towards the boundary of the interest region. Our method profits from the hierarchical nature of subdivision curves, since the unknowns of the optimization process are the few control points of the subdivision curve in the coarse representation and, at the same time, good approximations of the energies and their derivatives are obtained from the fine representation. We introduce a new region energy that guides the snake maximizing the contrast between the average intensity of the image within the snake and over the complement of the snake in a bounding box that does not change during the optimization. To illustrate the performance of our method we discuss the snakes associated with two classical subdivision schemes: the four point scheme and the cubic B-spline. Our experiments using synthetic and real images confirm that the proposed method is fast and robust.
\end{abstract}


\section{Introduction}
Active contours or snakes were introduced by Kass et al. in \cite{Kass87} as curves that slither within an image from some
initial position towards the contour of the object of interest. Snakes have become very popular in segmentation and tracking applications \cite{Bri07}, \cite{Del15} since they are very flexible and efficient.

The evolution of the snake is formulated as a minimization problem and the corresponding objective functions is usually known as snake energy. During the optimization process, the snake is iteratively updated
from a starting position until it reaches the minimum of the energy function. This energy measures the distance between the snake and the boundary of the object. It also controls some desirable properties of the final snake, such as the smoothness, the interpolation of distinguished points, etc. The quality of segmentation is determined by the choice of the energy terms and the starting position of the snake.

Kass et al. \cite{Kass87} originally formulated the snake energy as a linear combination of three terms: the {\em image energy}, which only depends on the image, the {\em internal energy}, which ensures the smoothness of the snake, and  the {\em constraint energy}, which allows that the user interacts with the snake. The specific definition of these energies depends on the application, on the nature of the image and also on the representation of the snake. The image energy guides the snake to the boundary of the interest object and it is the most important energy. It is usually defined as a weighted sum of a gradient based energy
\cite{Kass87}, \cite{Sta92}, that provides a good approximation of the contour of the object, and a region based energy \cite{Fig00}, \cite{The11}, that distinguishes  different homogeneous regions within the image.  Gradient based energies have a narrow zone of attraction in comparison with region based energies. Hence, the success of the segmentation depends on the
selection of the weight.

Snakes differ not only in the choice of the energy function but also in the representation of the curve. According to the representation, snakes may be classified as point snakes \cite{Kass87}, geodesic snakes \cite{Cas97}, \cite{Zha03}, \cite{App05}, \cite{Zha03} and  parametric snakes \cite{Fig00}, \cite{Bri00}, \cite{Shi11}, \cite{Del12a}.
Point-snakes simply consist of an ordered collection of points. This representation depends on a large number of parameters (the snake points), which makes the optimization expensive.  Geodesic snakes are
described as the zero level set of a higher-dimensional manifold. This type of active contours is very flexible topologically. In consequence, it is suitable for segmenting objects that have very variable shapes. A  drawback of geodesic snakes is that they are expensive from computational point of view. Parametric snakes are smooth curves written as a linear combination of a basis of functions.  The coefficients in this representation, known as control points, are few. This makes faster the optimization process. The downside of parametric snakes is that the parametrization restricts the shapes that can be approximated.

In this paper we focus on a particular class of parametric snakes: those generated from a subdivision scheme. Subdivision curves describe a contour by an initial discrete  and finite set of control points which, by the iterative application of refinement rules, becomes continuous in the limit.  Depending on the  particular choice of the subdivision mask, the continuous limit curve may have different degrees of smoothness. The main advantages of subdivision schemes are their simplicity of implementation,  the possibility to control their order of approximation, and their multiresolution property, which provides representations of the contour of a shape with varying resolutions.

\subsection{Related work}

The use of subdivision curves for segmentation was first proposed in \cite{Hug99}, where the so called {\it tamed snake} is introduced. This snake is generated by the classical four point subdivision scheme \cite{Dyn87}. The method incorporates image information considering the control points of the subdivision curve as mass points attracted by edges of the image. The four point subdivision scheme is also used in \cite{Shi11}  in combination with the gradient vector flow.  After every step of subdivision, the region energy of the subdivision polygon is reversely computed and a local adaptive compensation is carried out, in such a way that regions with high curvature are further subdivided, while flat regions remain unrefined.

In the context of image segmentation the most common snakes based on subdivision schemes are those producing B-spline type curves \cite{Bri00}, \cite{Fig00}, \cite{Jac04}, \cite{Del13a}. In \cite{Bri00} the snake is represented by  cubic B-spline basis functions.
The  initial B-spline is specified choosing node points instead of the B-spline control points in order to provide a more  intuitive user-interaction. To improve optimization speed and robustness a multiresolution approach is selected. This approach, based on an image pyramid, starts applying the optimization procedure at the coarsest level on a very small version of the image. After convergence, this solution is used as starting condition for  the next finer level.

In \cite{McI08} a segmentation method called {\em SketchSnakes} is proposed. The method combines a general subdivision curve snake with an initialization process  based on few sketch lines drawn by the user across the width of the target object. External image forces are computed at the points of the finer level curve  and then distributed, using weights derived from the original subdivision rules, among the points of the coarse level.  The positions of the control points are updated, new external forces are calculated and the process is repeated until an accurate solution is reached.

More recently, exponential B-spline have been introduced to construct snakes that reproduce circular and elliptical shapes \cite{Del12a}, \cite{Del12b}, \cite{Del13a}. In \cite{Bad17} subdivision snakes are obtained in a generic way using a multiscale approach to speed up the optimization process and improve robustness. Depending on the selected admissible subdivision mask, the snake may be interpolatory or reproduce trigonometric or polynomial curves. The multiscale approach facilitates to increase the number of points describing the curve as the algorithm progresses to the solution and, at each step, the scale of the image feature  is matched to the density of the sample of the curve.

\subsection{Our contribution}
The main contribution of this paper is a new region energy designed to maximize the contrast between the average intensity of the image within the snake and over the complement of the snake in a bounding box. This energy is simpler and computationally cheaper than other similar energies proposed in the literature \cite{Chan01},\cite{The11},\cite{Del12a}. In our region energy the  bounding box containing the object to be segmented does not change during the optimization. Moreover, the average intensity inside and outside the snake have neither to be estimated a priori nor to be included among the optimization parameters. Finally, in comparison with other methods, we are able to compute a better and more robust approximation of the region energy  using a method  to obtain a pixel-level discretization of the snake.  Our method produces good approximations of the region energy for images of either low or high resolution.

%

\section{Subdivision curves}
\label{sec:subdivision}

\subsection{Linear uniform stationary subdivision schemes}
Denote by $\bP^0$ a polygon in the plane. A {\it subdivision scheme} is a procedure that refines $\bP^0$  producing a sequence of polygons $\bP^1, \bP^2, \ldots$ with an increasing number of vertices. A {\it linear, stationary, uniform} and {\it binary} subdivision scheme is based on the application of a refinement rule
\begin{equation}
\bp_i^{k+1} = \sum_{j \in \ZZ} a_{i-2j} \bp_j^k
\label{eq:binsub}
\end{equation}
that computes the vertices ${\mathbf p}_i^k=(x_i^k,y_i^k)$  of the polygon $\bP^k$ in the step $k$ as a linear function of the vertices of $\bP^{k-1}$. The coefficients $\mathbf{a} = \{a_i \in \RR, i\in \ZZ\}$ in \eqref{eq:binsub}) are a called {\it subdivision mask}. In practice, only a finite number of coefficients are different from zero. The subdivision scheme converges if the sequence of piecewise linear functions $\mathbf{f}^k(t)$ which satisfies the interpolation conditions
\begin{equation}
\label{eq:piecewise_seq}
\mathbf{f}^k\left(\frac{i}{2^k}\right) = {\bp}_i^k, \quad i \in \ZZ
\end{equation}
{\it converges uniformly}. Denote by $\br(t)=(x(t),y(t))$ the continuous limit function
\begin{equation} \label{eq:ra}
\br(t)=\lim_{k\rightarrow \infty} \mathbf{f}^k(t).
\end{equation}
This limit exists as long as the subdivision scheme applied to the functional data $\mathbf{\delta}=\{\delta_{i,0},i \in \mathbb{Z}\}=\{\ldots,0,0,1,0,0,\ldots\}$ converges. The corresponding limit function $\varphi(t)$ is called {\it basic limit function} and satisfies the {\it refinement equation}
\begin{equation}
\varphi(t)=\sum_{j \in \mathbb{Z}} a_j \varphi(2t-j).
\label{eq:refi}
\end{equation}
Due to the linearity of the subdivision rules, the subdivision curve $\br(t)$  can be written as a linear combination of the integer shifts of $\varphi(t)$,
\begin{equation}
\br(t) = \sum_{j \in \mathbb{Z}} \bp_j^0 \varphi(t-j).
\label{eq:ra_basicdef}
\end{equation}
This representation can be used to define also the tangent vector to the curve (as well as the normal vector),
\begin{equation}
\frac{\mathrm{d}}{\mathrm{d}t}\br(t) = \sum_{j \in \mathbb{Z}} \bp_j^0 \frac{\mathrm{d}}{\mathrm{d}t}\varphi(t-j).
\label{eq:tanra_basicdef}
\end{equation}

Even more, because of \eqref{eq:binsub}, for any $k \geq 0$ the subdivision curve may be expressed as,
\begin{equation}
\br(t)= \sum_{j \in \mathbb{Z}} \bp_j^k \varphi(2^kt-j).
\label{eq:rak}
\end{equation}

From \eqref{eq:rak} we observe that,
\begin{equation}
\br\left(\frac{i}{2^k}\right) = \sum_{j \in \mathbb{Z}} \bp_j^0 \varphi\left(\frac{i}{2^k} - j\right) = \sum_{j \in \mathbb{Z}} \bp_j^k \varphi(i-j).
\label{eq:evalr}
\end{equation}

In case of the interpolatory subdivision schemes, as $\varphi(i-j) = \delta_j^i$, then it holds that,
\begin{equation}
\br\left(\frac{i}{2^k}\right) = \sum_{j \in \mathbb{Z}} \bp_j^k \delta_j^i = \bp_i^k.
\label{evalr_int}
\end{equation}

Also, because of \eqref{eq:piecewise_seq} and \eqref{eq:ra} it holds,
\begin{equation} \label{eq:limit_pos}
 \br\left(\frac{i}{2^k}\right) = \lim_{q\rightarrow\infty} \mathbf{f}^{k+q}\left(\frac{2^{q}i}{2^{k+q}}\right) = \lim_{q\rightarrow\infty} \bp_{2^q i}^{k+q}, \qquad q\in\NN.
\end{equation}
In case of interpolatory schemes the center and right hand expressions are limits of constant sequences that converge to the point $\bp_i^k$ as stated in \eqref{evalr_int}. Otherwise, the right hand sequence allows to analyze  the exact value of $\br\left(i/2^k\right)$, as it is done in \eqref{eq:limit_pos_cBS}.

\begin{remark} \label{rem:evalc_vs_subdiv}
A subdivision curve $\br(t)$ is usually represented by a polygon $\bP^k$ whose vertices are obtained after some refinements of an initial polygon $\bP^0$. As $k$ increases the polygon $\bP^k$ provides a better approximation of $\br(t)$. In this work we approximate $\br(t)$ by the polygon $\{\br(i/2^k)\}$ whose vertices are on the curve. For interpolatory subdivision schemes both $\bP^k$ and $\{\br(i/2^k)\}$ polygons are the same, but they are different in the case of non-interpolatory schemes.
\end{remark}

The expression \eqref{eq:evalr} allows two options for the evaluation of the curve in dyadic parameters. The first one uses the initial control polygon $\bP^0$ and computes the evaluations $\varphi\left(\frac{i}{2^k} - j\right)$. The second subdivides the initial polygon $k$ times and uses the polygon $\bP^k$. We choose the first option, since the evaluation of $\varphi\left(\frac{i}{2^k} - j\right)$ for any $i,j \in \ZZ$ and $k \in \NN$ can be done {\em a priori} storing the results in a lookup table. The evaluation of the basic function at dyadic parametric values $\varphi\left(\frac{i}{2^k} - j\right)$ for any $i,j \in \ZZ$ and $k \in \NN$ can be computed as the subdivision of the polygon with vertices ${\mathbf P}^0 = \left\{ (i,\delta_i^0), \; i \in \ZZ \right\}$ by $k$ times (see Figures \ref{fig:basic_4pt} and \ref{fig:basic_BS}).

In this work we are interested in {\it closed curves} that approximate the boundary of a region in a digital image. Hence, if the curve $\br(t)$  has  a control polygon with $M$ vertices $\bp_i^0,\; i=0,\ldots,M-1$, then the polygon is periodically extended assuming that $\bp_{i+M}^0=\bp_{i}^0$ for all $i \in \mathbb{Z}$. Under this assumption, we obtain from \eqref{eq:ra_basicdef})
\begin{align}
\br(t+M) &= \sum_{j \in \mathbb{Z}} \bp_j^0 \varphi(t+M-j)=\sum_{i \in \mathbb{Z}} \bp_{i+M}^0 \varphi(t-i) \\ \notag
 &= \sum_{i \in \mathbb{Z}} \bp_{i}^0 \varphi(t-i)
=\br(t),
\end{align}
i.e., the subdivision curve is periodic with period $M$. For more details about subdivision schemes, see \cite{Dyn02}.

\subsection{Exact evaluation of linear uniform stationary subdivision schemes}

To illustrate the performance of our method we discuss in details two classical subdivision schemes: the four point and the cubic B-spline.  The first one is interpolatory while the second is non-interpolatory.  Since the interpolation of the points provided by the user is the more natural starting point for the snake, we explain in the case of the cubic B-spline how to compute the initial control polygon in such a way that the B-spline  curve passes through the given set of points.

In the next sections we describe the subdivision schemes chosen to generate the snake curves. Moreover, we provide the expressions needed to evaluate the subdivision curves and their derivatives at dyadic parameters.

\subsubsection{Four points subdivision scheme}

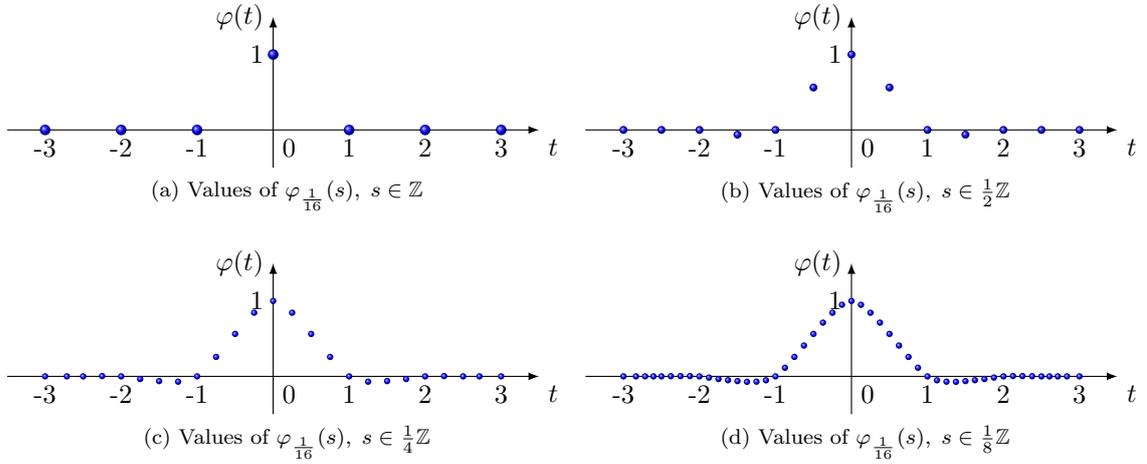
\begin{figure}[h]
  \centering
\subfloat[Values of $\varphi_{\frac{1}{16}}(s), \; s \in \ZZ$]{
  \begin{tikzpicture}
  \draw[tip] (0,-.5) -- (0,1.5) node[left] {$\varphi(t)$};
  \draw[tip] (-3.5,0) -- (3.5,0) node[below right] {$t$};
  \node[left] at (0,1) {1};
  \foreach \i in {-3,-2,-1,1,2,3} \node[below] at (\i,0) {\i};
  \node[below right] at (0,0) {0};

  \foreach \x/\y [count=\i] in {-3/0,-2/0,-1/0, 0/1, 1/0,2/0,3/0}
  \coordinate (p\i) at (\x,\y);

  \foreach \j in {1,...,7}
  \shade[ball color=blue] (p\j) circle (2pt);

  \end{tikzpicture}
  }
  \subfloat[Values of $\varphi_{\frac{1}{16}}(s), \; s \in \frac12\ZZ$]{
    \begin{tikzpicture}
    \draw[tip] (0,-.5) -- (0,1.5) node[left] {$\varphi(t)$};
    \draw[tip] (-3.5,0) -- (3.5,0) node[below right] {$t$};
    \node[left] at (0,1) {1};
    \foreach \i in {-3,-2,-1,1,2,3} \node[below] at (\i,0) {\i};
    \node[below right] at (0,0) {0};

  \foreach \x/\y [count=\i] in {-4/0,-3/0,-2/0,-1/0, 0/1, 1/0,2/0,3/0,4/0}
    \coordinate (p\i) at (\x,\y);

\ssfourpoints(9,1)
    \foreach \j in {3,...,15}
    \shade[ball color=blue] (p\j) circle (1.5pt);

    \end{tikzpicture}
    }

    \subfloat[Values of $\varphi_{\frac{1}{16}}(s), \; s\in \frac14\ZZ$]{
      \begin{tikzpicture}
      \draw[tip] (0,-.5) -- (0,1.5) node[left] {$\varphi(t)$};
      \draw[tip] (-3.5,0) -- (3.5,0) node[below right] {$t$};
      \node[left] at (0,1) {1};
      \foreach \i in {-3,-2,-1,1,2,3} \node[below] at (\i,0) {\i};
      \node[below right] at (0,0) {0};

  \foreach \x/\y [count=\i] in {-4/0,-3/0,-2/0,-1/0, 0/1, 1/0,2/0,3/0,4/0}
      \coordinate (p\i) at (\x,\y);

\ssfourpoints(9,2)
      \foreach \j in {5,...,29}
      \shade[ball color=blue] (p\j) circle (1.2pt);

      \end{tikzpicture}
      }
      \subfloat[Values of $\varphi_{\frac{1}{16}}(s), \; s \in \frac18\ZZ$]{
        \begin{tikzpicture}
        \draw[tip] (0,-.5) -- (0,1.5) node[left] {$\varphi(t)$};
        \draw[tip] (-3.5,0) -- (3.5,0) node[below right] {$t$};
        \node[left] at (0,1) {1};
        \foreach \i in {-3,-2,-1,1,2,3} \node[below] at (\i,0) {\i};
        \node[below right] at (0,0) {0};

  \foreach \x/\y [count=\i] in {-4/0,-3/0,-2/0,-1/0, 0/1, 1/0,2/0,3/0,4/0}
        \coordinate (p\i) at (\x,\y);

        \ssfourpoints(9,3)
        \foreach \j in {9,...,57}
        \shade[ball color=blue] (p\j) circle (1.2pt);

        \end{tikzpicture}
        }
    \caption{Generating the values of the basic function for the 4-point subdivision scheme.}
\label{fig:basic_4pt}
\end{figure}

The four point subdivision scheme \cite{Dyn87}, also known as DLG, is a linear, stationary and uniform subdivision scheme, depending on a tension parameter $\omega$. The rules that define this scheme are
\begin{align}
\bp_{2i}^{k+1} &= \bp_i^k \label{4pt_evenrule}\\
\bp_{2i+1}^{k+1} &= \left(\omega + \frac{1}{2}\right)\left(\bp_i^k + \bp_{i+1}^k\right) - \omega \left(\bp_{i-1}^k + \bp_{i+2}^k\right)
\label{4pt_oddrule}
\end{align}
The scheme is {\it interpolatory} since the rule \eqref{4pt_evenrule} implies that the set of points of the step $k+1$ contains the points of the previous step. Hence, the control points $\bp_i^0$ are interpolated by functions $\mathbf{f}^k(t)$ for all $k \geq 0$ and thus they belong to the limit curve that we denote by  $\br_{\omega}(t)$, to recall that it depends on the free parameter $\omega$. This curve is continuous if $\omega$ is in the interval $(0,\frac{1}{4})$ and it has a continuous tangent vector when $\omega \in \left(0,\frac{\sqrt{5} - 1}{8}\right)$. The basic limit function $\varphi_{\omega}$ of the DLG scheme has support $[-3,3]$, as it is shown in Figure \ref{fig:basic_4pt} for $\omega = \frac{1}{16}$.

Since the subdivision rules \eqref{4pt_evenrule}-\eqref{4pt_oddrule} depend on 4 points to get a closed curve we need a closed polygon composed by $M+3$ points ${\mathbf P}^0=\{\bp_{-1}^0, \bp_0^0, \ldots, \bp_{M-1}^0, \bp_{M}^0, \bp_{M+1}^0\}$,
where  $\bp_{-1}^0 = \bp_{M-1}^0$, $\bp_M^0=\bp_0^0$ and $\bp_{M+1}^0 = \bp_1^0$. Hence, if the initial polygon has $M$ vertices, then expression \eqref{eq:ra_basicdef} for the subdivision curve is reduced to
\begin{equation}
\br_{\omega}(t)= \sum_{j=-1}^{M+1} \bp_j^0 \varphi_{\omega}(t-j),\quad 0 \leq t \leq M
\label{eq:4pt_raredu}
\end{equation}
where $\varphi_{\omega}(t)$ is the basic limit function of the four point subdivision scheme with parameter $\omega$.

\medskip

If we denote by $\mathbf{t}_i^k$ the tangent vector to the subdivision curve at $\bp_i^k$, then it holds:
\begin{equation}
\mathbf{t}_i^k = \left.\frac{\mathrm{d}}{\mathrm{d}t}\br_{\omega}(t)\right\rvert_{t=\frac{i}{2^k}}
= \frac{2^k}{1-4\omega}\left(\frac{1}{2}(\bp_{i+1}^k-\bp_{i-1}^k)- \omega(\bp_{i+2}^k-\bp_{i-2}^k)\right).
\label{eq:4pt_tanik}
\end{equation}

\begin{remark} \label{rem:tan}
In practice, we use \eqref{eq:4pt_tanik} only to compute $\varphi'(i/2^k)$, obtained with $\{\bp_i^0 = (i, \delta_i^0)\}$. The tangent vector for any subdivision curve is computed by using \eqref{eq:tanra_basicdef}.
\end{remark}

Recall that for $i=0,\ldots,2^kM-1$ the subdivision curve $\br_{\omega}(t)$ of the DLG scheme satisfies
\begin{equation}
\br_{\omega}\left(\frac{i}{2^k}\right)=\bp_i^k, \qquad
\frac{d \br_{\omega}}{dt}\left(\frac{i}{2^k}\right) = \mathbf{t}_i^k.
\label{eq:4pt_piktik}
\end{equation}
In what follows we use the notation $\mathbf{t}_i^k = (tx_i^k,ty_i^k)$.

The best value for the parameter $\omega$ with respect to the regularity of the limit curve is $\omega = \tfrac{1}{16}$ \cite{Dyn:LaNL_SS}. In the rest of this paper we consider only this case.


\subsubsection{Cubic B-spline subdivision scheme}

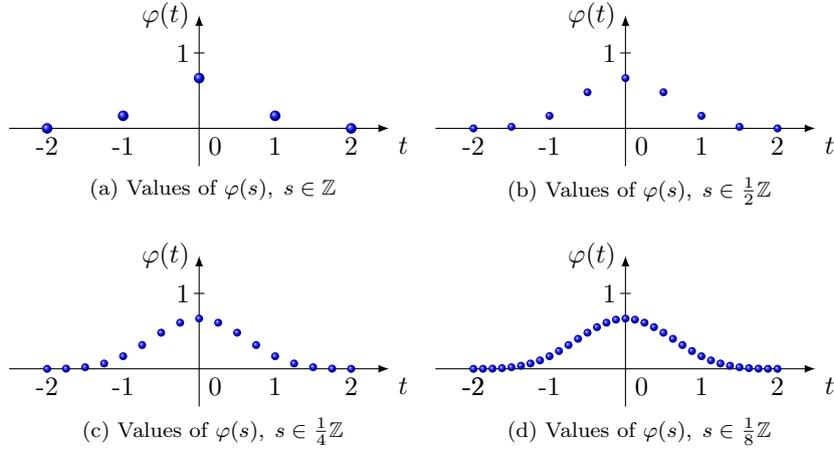
\begin{figure}[h]
  \centering
\subfloat[Values of $\varphi(s), \; s \in \ZZ$]{
  \begin{tikzpicture}
  \draw[tip] (0,-.5) -- (0,1.5) node[left] {$\varphi(t)$};
  \draw[tip] (-2.5,0) -- (2.5,0) node[below right] {$t$};
 \draw (0,1) node[left] {1} +(-.07,0) -- +(.07,0);
  \foreach \i in {-2,-1,1,2} \node[below] at (\i,0) {\i};
  \node[below right] at (0,0) {0};

  \foreach \x/\y [count=\i] in {-3/0,-2/0,-1/0, 0/1, 1/0,2/0,3/0}
  \coordinate (p\i) at (\x,\y);

\foreach \h in {2,...,6} {%
  \pgfmathtruncatemacro{\p}{\h-1}
  \pgfmathtruncatemacro{\j}{\h+1}
  \coordinate (n\h) at (barycentric cs:p\p=1/6,p\h=2/3,p\j=1/6);
}
\foreach \j in {2,...,6} \shade[ball color=blue] (n\j) circle (2pt);
  \end{tikzpicture}
  }
  \subfloat[Values of $\varphi(s), \; s \in \frac12\ZZ$]{
    \begin{tikzpicture}
    \draw[tip] (0,-.5) -- (0,1.5) node[left] {$\varphi(t)$};
    \draw[tip] (-2.5,0) -- (2.5,0) node[below right] {$t$};
    \draw (0,1) node[left] {1} +(-.07,0) -- +(.07,0);
    \foreach \i in {-2,-1,1,2} \node[below] at (\i,0) {\i};
    \node[below right] at (0,0) {0};

  \foreach \x/\y [count=\i] in {-3/0,-2/0,-1/0, 0/1, 1/0,2/0,3/0}
    \coordinate (p\i) at (\x,\y);

\ssBS(7,1)
\foreach \h in {3,...,11} {%
  \pgfmathtruncatemacro{\p}{\h-1}
  \pgfmathtruncatemacro{\j}{\h+1}
  \coordinate (n\h) at (barycentric cs:p\p=1/6,p\h=2/3,p\j=1/6);
}
\foreach \j in {3,...,11} \shade[ball color = blue] (n\j) circle (1.5pt);
    \end{tikzpicture}
    }

    \subfloat[Values of $\varphi(s), \; s\in \frac14\ZZ$]{
      \begin{tikzpicture}
      \draw[tip] (0,-.5) -- (0,1.5) node[left] {$\varphi(t)$};
      \draw[tip] (-2.5,0) -- (2.5,0) node[below right] {$t$};
      \draw (0,1) node[left] {1} +(-.07,0) -- +(.07,0);
      \foreach \i in {-2,-1,1,2} \node[below] at (\i,0) {\i};
      \node[below right] at (0,0) {0};

  \foreach \x/\y [count=\i] in {-3/0,-2/0,-1/0, 0/1, 1/0,2/0,3/0}
      \coordinate (p\i) at (\x,\y);

\ssBS(7,2)
\foreach \h in {5,...,21} {%
  \pgfmathtruncatemacro{\p}{\h-1}
  \pgfmathtruncatemacro{\j}{\h+1}
  \coordinate (n\h) at (barycentric cs:p\p=1/6,p\h=2/3,p\j=1/6);
}
\foreach \j in {5,...,21} \shade[ball color=blue] (n\j) circle (1.5pt);
      \end{tikzpicture}
      }
      \subfloat[Values of $\varphi(s), \; s \in \frac18\ZZ$]{
        \begin{tikzpicture}
        \draw[tip] (0,-.5) -- (0,1.5) node[left] {$\varphi(t)$};
        \draw[tip] (-2.5,0) -- (2.5,0) node[below right] {$t$};
        \draw (0,1) node[left] {1} +(-.07,0) -- +(.07,0);
        \foreach \i in {-2,-2,-1,1,2} \node[below] at (\i,0) {\i};
        \node[below right] at (0,0) {0};

  \foreach \x/\y [count=\i] in {-3/0,-2/0,-1/0, 0/1, 1/0,2/0,3/0}
        \coordinate (p\i) at (\x,\y);

        \ssBS(7,3)
        \foreach \h in {9,...,41} {%
          \pgfmathtruncatemacro{\p}{\h-1}
          \pgfmathtruncatemacro{\j}{\h+1}
          \coordinate (n\h) at (barycentric cs:p\p=1/6,p\h=2/3,p\j=1/6);
        }
\foreach \j in {9,...,41} \shade[ball color=blue] (n\j) circle (1.5pt);
        \end{tikzpicture}
        }

  \caption{Generating the values of the basic function for the cubic B-spline subdivision scheme.}
\label{fig:basic_BS}
\end{figure}

This linear, stationary and uniform subdivision scheme is defined by the rules:
\begin{align}
\bp_{2i}^{k+1} &= \tfrac{1}{8}\bp_{i-1}^k + \tfrac{6}{8}\bp_{i}^k + \tfrac{1}{8}\bp_{i+1}^k \label{BS_evenrule}\\
\bp_{2i+1}^{k+1} &= \tfrac{1}{2}\bp_i^k + \tfrac{1}{2}\bp_{i+1}^k
\label{BS_oddrule}
\end{align}
and generates as limit a cubic B-spline curve that is $C^2$-continuous. The basic limit function $\varphi$ for this scheme has support $[-2,2]$, as it is shown in Figure \ref{fig:basic_BS}.

Since this scheme is not interpolatory, the points in ${\mathbf P}^k$ don't belong to the limit curve. Nevertheless, following \eqref{eq:limit_pos} it can be proved \cite{deRose} that:
\begin{equation}  \label{eq:limit_pos_cBS}
 \br\left(\frac{i}{2^k}\right) = \tfrac{1}{6}\bp_{i-1}^{k} + \tfrac{4}{6}\bp_{i}^{k} + \tfrac{1}{6}\bp_{i+1}^{k}.
\end{equation}

\begin{remark}
Recall that in this work, we don't compute the values $\br\left(i/2^k\right)$ using \eqref{eq:limit_pos_cBS} for any polygon $\bP^0$ and its refinements. Instead, we store in a lookup table the pre-computed values of $\varphi\left(\frac{i}{2^k}\right)$ for $i \in \ZZ$ and a previous fixed value of $k \in \NN$, obtained from the initial data  ${\bP}^0 = \left\{ (i, \delta_i^0), \; i \in \ZZ \right\}$. Then, we use:
\[
\br\left(\frac{i}{2^k}\right) = \sum_{j \in \mathbb{Z}} \bp_j^0 \varphi\left(\frac{i}{2^k} - j\right).
\]
Taking into account that $\bP^0$ changes during the optimization of the snake the previous strategy reduces the computational cost.
\end{remark}

Since the subdivision rules \eqref{BS_evenrule}-\eqref{BS_oddrule} depend on 3 points to get a closed curve we need a closed polygon composed by $M+2$ points ${\mathbf P}^0=\{\bp_{-1}^0, \bp_0^0, \ldots, \bp_{M-1}^0, \bp_{M}^0\}$, where  $\bp_{-1}^0 = \bp_{M-1}^0$ and $\bp_M^0=\bp_0^0$. Hence, if the initial polygon has $M$ vertices, then the expression \eqref{eq:ra_basicdef} for the subdivision curve is reduced to,
\begin{equation}
\br(t)= \sum_{j=-1}^{M} \bp_j^0 \varphi(t-j),\quad 0 \leq t \leq M
\label{eq:BS_raredu}
\end{equation}
where $\varphi(t)$ is the basic limit function of cubic B-spline subdivision scheme.


If we denote by $\mathbf{t}_i^k$ the tangent vector to the subdivision curve at $\bp_i^k$, then,
\begin{equation}
\mathbf{t}_i^k = \left.\frac{\mathrm{d}}{\mathrm{d}t}\br(t)\right\rvert_{t=\frac{i}{2^k}}
	= \frac{1}{2}\left(\bp_{i+1}^k - \bp_{i-1}^k\right).
\label{eq:BS_tanik}
\end{equation}

The use of this expression follows the same argument as in {\em Remark \ref{rem:tan}}.

\medskip

Approximating subdivision curves do not interpolate their control points. Since control points are the degrees of freedom, this property makes approximating subdivision snakes less intuitive in interactive segmentation than interpolating subdivision snakes. To overcome this limitation, we explain now how to compute the control points $\widetilde{\bP}^0$ of the cubic B-spline curve interpolating the vertices of the initial polygon introduced by the user. This new strategy, different from other works as \cite{Bri00}, unifies the treatment of approximating and interpolating subdivision snakes, making the initialization and interaction
more user-friendly.

\begin{figure}[h]
  \centering
  \begin{tikzpicture}[scale=.7]
  \foreach \x/\y [count = \i] in {-.75/-.75, 3.75/-.75, 3.75/3.75, -.75/3.75}
  \coordinate (p\i) at (\x,\y);
\node[left] at (p1) {$\mathbf{\widetilde{p}}_0^0$};
\node[right] at (p2) {$\mathbf{\widetilde{p}}_1^0$};
\node[right] at (p3) {$\mathbf{\widetilde{p}}_2^0$};
\node[left] at (p4) {$\mathbf{\widetilde{p}}_3^0$};

\draw[dotted] (p1) -- (p2) -- (p3) -- (p4) -- cycle;
\foreach \i in {1,...,4}			
\shade[ball color = red] (p\i) circle (2pt);
\ssBS(4,5)
\draw (p1) \foreach \j in {2,...,\i}{ -- (p\j)} -- cycle;
			
  \foreach \x/\y [count = \i] in {0/0, 3/0, 3/3, 0/3}
  \coordinate (ip\i) at (\x,\y);
\node[above right] at (ip1) {$\bp_0^0$};
\node[above left] at (ip2) {$\bp_1^0$};
\node[below left] at (ip3) {$\bp_2^0$};
\node[below right] at (ip4) {$\bp_3^0$};

\foreach \i in {1,...,4}			
\shade[ball color = blue] (ip\i) circle (2pt);
  \end{tikzpicture}
  \caption{Interpolation of given set of points $\bP^0$ by the cubic B-spline subdivision scheme}
  \label{fig:BS_interpolatory}
\end{figure}
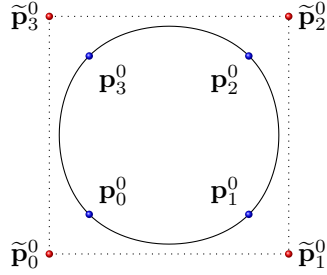

The control polygon $\widetilde{\bP}^0$ depends linearly on the polygon $\bP^0$ as it is shown in the following.

\begin{theorem}
The subdivision curve generated by the cubic B-spline scheme that interpolates the set of points $\bP^0$ has control points $\widetilde{\bP}^0$ given by,
\begin{equation}
  \label{eq:pseudo_int}
 \widetilde{\bP}^0 = \mathcal{A} {\bP}^0,
\end{equation}
where the elements of matrix $\mathcal{A} = [a_{s,t}]_{s=0,t=0}^{M-1,M-1}$ are given by:
\begin{equation*}
a_{s,t} =  \begin{cases}
\frac1M + \frac3M \cos(s\pi) +  \frac2M \displaystyle \sum_{j=1}^{\frac{M}{2} -1} \left( \frac23 + \frac13 \cos(2j\pi/M) \right)^{-1} \cos(2(s-t)\pi/M),  & \text{for} \quad M\mod{2} = 0, \\
\frac1M + \frac2M \displaystyle \sum_{j=1}^{\lfloor\frac{M}{2}\rfloor} \left( \frac23 + \frac13 \cos(2j\pi/M) \right)^{-1} \cos(2(s-t)\pi/M), & \text{for} \quad M\mod{2} = 1.
\end{cases}
\end{equation*}

\end{theorem}

See the proof in {\em Appendix \ref{sec:interp_op}}. Figure \ref{fig:BS_interpolatory}  shows the points $\bP^0$ to be interpolated by the cubic B-spline and the control polygon $\widetilde{\bP}^0$ computed using \eqref{eq:pseudo_int}.

\section{Snake Energies}
\label{sec:energies}

In the literature, the evolution of the snake is driven by the minimization of several energies that measure the proximity between the snake and the boundary $\partial\Gamma$ of a bounded region $\Gamma$ in a digital image and also some desirable properties of the final curve like the smoothness, the interpolation of distinguished points and so on.

Since our snake is a subdivision curve, the total energy $E_{\text{snake}}$, depends on the initial control polygon $\bP^0$. The control polygon $\bP_{*}^0$ of the optimal snake is computed as:
\begin{equation}
\bP_{*}^0= arg\;\min_{\bP^0}\;E_{\text{snake}}(\bP^0).
\label{optpar}
\end{equation}

In this paper we assume that the region of interest $\Gamma$ to be segmented is dark in comparison to the background. Hence, the energy functionals related with the image are designed to detect dark objects on a brighter background.
All the energies are defined by integrals of functions which are computed approximately. To obtain good approximations we use a large sample of points on the subdivision curve. In the following sections we develop the expressions for each energy.

\begin{remark} \label{rem:coord_syst}
Images are represented in a system of coordinates defined by rows and columns, like the indexing of a matrix. Thus, if a pixel has coordinates $(x,y)$, the $x$-coordinate refers to the row and the $y$-coordinate refers to the column (see for example Figure \ref{fig:pixels_reg_energy}). This does not affect the definition and use of the subdivision schemes, since each coordinate in \eqref{eq:ra_basicdef} works independently.
\end{remark}

\subsection{Gradient energy}

If $I(x,y)$ denotes the image intensity at a pixel with coordinates $(x,y)$ and $\br(t)=(x(t),y(t))$ is a parametric curve living on the image  for $t \in[0,M]$, the simplest image energy is the {\it gradient magnitude energy} $E_{mag}$ given by:
\begin{equation}
E_{mag}(\br(t))= - \int_{0}^M \|\nabla I (\br(t))\|^2\;dt
\label{Emag}
\end{equation}
where $\lVert\nabla I (\br(t))\rVert^2 = \left(\frac{\partial I}{\partial x}(x(t),y(t))\right)^2 \; + \; \left(\frac{\partial I}{\partial y}(x(t),y(t))\right)^2$.
Since the gradient magnitude energy only depends on the magnitude of the gradient vector, the minimization of $\eqref{Emag}$ can misguide the snake to a neighboring object if the initial approximation is not very close to the boundary of interest. To overcome this limitation, several alternatives energies has been proposed, like balloon forces \cite{Coh91}, gradient vector-fields  \cite{Xu98}, \cite{Jac01},\cite{Jac04} or  multiresolution approaches \cite{Bri00}.

In this paper we use the gradient based image energy $E_{grad}$ proposed in \cite{Jac01}. The idea behind this approach is the following. If we travel around the ground truth boundary curve $\partial \Gamma$ in counterclockwise direction, then $\Gamma$ is always on the ``left", i.e in the direction of $-\nabla I$. Hence, we pull the snake in the direction of $\partial \Gamma$, requiring  the normal to snake at any point to be parallel to $-\nabla I$ at the same point. More precisely, if we denote by $\mathbf{n}(t)$ the inward unit normal to snake  at the point $\br(t)$, then the new energy $E_{grad}$, which takes into account not only the magnitude of the image gradient but also its direction is given by:
\begin{equation}
E_{grad}(\br(t))= - \int_{0}^M \langle\nabla I (\br(t)) , \left\lVert  \frac{\mathrm{d}\br(t)}{\dt} \right\rVert \mathbf{n}(t)\rangle \;dt
\label{Egrad}
\end{equation}
where $\langle \cdot,\cdot \rangle$ is the usual scalar product  and $\frac{\mathrm{d}\br(t)}{\dt}$ denotes the tangent to $\br(t)$. Expanding \eqref{Egrad} we obtain
\begin{equation}
E_{grad}(\br(t)) = - \int_{0}^M \left(\frac{\partial I}{\partial x}(x(t),y(t))\frac{\dy(t)}{\dt} - \frac{\partial I}{\partial y}(x(t),y(t))\frac{\dx(t)}{dt}\right) \;dt. \label{Egrad2}
\end{equation}

To obtain good approximations of the energies (and their derivatives with respect to the coordinates of control points) we use a large sample of points on the subdivision curve. More precisely, given the initial polygon ${\mathbf P}^0=\{\bp_0^0,\ldots,\bp_{M-1}^0\}$, we select $k$ (in our experiments we take $k=4$ or $k=5$) and we use \eqref{eq:ra_basicdef} to generate $2^kM$ points $\br(i/2^k)$, $i=0,\ldots,2^kM-1$ on the subdivision curve. Moreover, we apply bilinear interpolation on the gradient of the image to compute $\nabla I(\br(i/2^k))$. Finally, we approximate the energy  substituting the integral in \eqref{Egrad2} by the average of values of the integrand over the sample of $2^kM$ points on the subdivision curve corresponding to parameter values $\frac{i}{2^k}$, $i=0,\ldots,2^kM-1$.

Taking into account \eqref{eq:ra_basicdef} we obtain\footnote{Recall that the indices of the inner summations depend on the choice of the subdivision scheme, see \eqref{eq:4pt_raredu} and \eqref{eq:BS_raredu}.} the following approximation of \eqref{Egrad2}:
\begin{equation}
E_{grad}({\mathbf P}^0)\approx \frac{1}{2^kM} \sum_{i=0}^{2^kM-1}\;\left[ \frac{\partial I}{\partial y}\left(\sum_{j\in\ZZ} \bp_j^0 \varphi(\tfrac{i}{2^k}-j)\right)tx_i^k - \frac{\partial I}{\partial x}\left(\sum_{j\in\ZZ} \bp_j^0 \varphi(\tfrac{i}{2^k}-j)\right)ty_i^k\right]
\label{Egrad3}
\end{equation}
where $tx_i^k = \frac{\dx}{\dt}\left(i/2^k\right)$ and $ty_i^k = \frac{\dy}{\dt}\left(i/2^k\right)$, so that $\left(tx_i^k,ty_i^k\right) = \frac{\mathrm{d}\br}{\dt}\left(i/2^k\right)$.

It should be noticed that the right hand side of \eqref{Egrad3} is a function of the coordinates of the initial control points ${\mathbf P}^0$.

\subsection{Region energy}
The main limitation of gradient based energy \eqref{Egrad} is that its zone of attraction is limited, since the gradient is small as long as we move away from $\partial \Gamma$. To face this problem several region energies have been introduced in the literature \cite{Sta92}, \cite{Ches99}, \cite{Chan01},\cite{The11}, \cite{Del12a}. Some of them use statistical information to identify different regions \cite{Jac01},\cite{Jac04},\cite{The11}.
Inspired in the energies proposed in  \cite {Del12a} and \cite{The11} we introduce in this work a simple region energy $E_{reg}$, designed to maximize the contrast between the average intensity of the pixels within the snake and the average intensity in the region outside the snake and inside a given bounding box.

Assuming that $\Omega$, the region enclosed by the snake $(x(t),y(t)),\,t\in [0,M]$, is contained in a rectangular region $R$, we denote by $\left| R \right|$ the area of $R$ (which is a constant) and by  $\left| \Omega \right|$ the area of $\Omega$ (which may vary, while the snake evolves). The new region energy, $E_{reg}$, to be minimized is,
\begin{equation}
E_{reg}(\bP^0) := -\left({\frac {\int\int_{\Omega} I(x,y) \dx\dy  }{\lvert\Omega\rvert }}
-{\frac {\int\int_{R\setminus \Omega} I(x,y) \dx\dy }{\lvert R \rvert - \lvert \Omega \rvert }} \right)^2.
\label{Ereg}
\end{equation}

Observe that minimizing $E_{reg}$ is equivalent to  maximize the difference between the average intensity inside $\Omega$  and the average intensity in the complement of $\Omega$ in $R$.

Let us introduce the following notation,
\[
I_{\Omega}:=\int\int_{\Omega} I(x,y) \dx\dy  \qquad
\text{and} \qquad
I_{R}:=\int\int_{R} I(x,y) \dx\dy.
\]
Then, the region energy may be written as,

\begin{equation}
E_{reg}(\bP^0) = -\left(\frac{I_{\Omega}}{\lvert \Omega \rvert} - \frac{I_R-I_{\Omega}}{\lvert R \lvert - \lvert \Omega \rvert} \right)^2.
\label{Ereg1}
\end{equation}

Since region energies are usually expressed as integrals of a function over the domain $\Omega$ enclosed by the snake, some authors propose the use of Green's theorem to rewrite the 2D integrals as a line integral along the snake \cite{Del12a},\cite{Del13a},\cite{Jac04}. In particular, if we apply it to the function $I(x,y)$ we obtain
\begin{equation}
I_\Omega = \int\int_{\Omega} I(x,y) \dx\dy = \int_{\partial \Omega} I_1(x,y) \dy = -\int_{\partial \Omega} I_2(x,y) \dx
\label{Green}
\end{equation}
where
\begin{equation}
 I_1(x,y) = \int_{-\infty}^x I(\tau,y) \mathrm{d}\tau \quad \text{and} \quad
 I_2(x,y) = \int_{-\infty}^y I(x,\tau) \mathrm{d}\tau.
 \label{I1I2}
\end{equation}
Thus, if $\partial \Omega$ is parametrized by $\br(t)=(x(t),y(t)),\;0 \leq t \leq M$, then from \eqref{Green}) it holds
\begin{equation}
I_\Omega = \int_{0}^M I_1(x(t),y(t))\; \frac{\dy(t)}{\dt} \dt\\
                                = -\int_{0}^M I_2(x(t),y(t))\;\frac{\dx(t)}{\dt} \dt.
\label{intI}
\end{equation}

This approach reduces significantly the computational cost, but in our experiments we have found that large errors may be introduced when we use it to compute the integrals in \eqref{Ereg}, in a digital images context.

In those works, the line integrals \eqref{I1I2} are approximated using a sample of points on the snake and summing up the contributions of column or row image pixel strips corresponding to each point on the snake.  But even if the snake is parametrized by a multiple of the arc length, the distribution on the image of the sample of points may be very irregular. For instance, if the image has low resolution then some points may belong to the same pixel overestimating the value of the integral. On the contrary, if the  image has high resolution then those rows or columns of $\Omega$ without any point of the sample do not contribute to the computation producing an underestimate of the integral. We propose instead a sort of {\em rasterization} of $\partial \Omega$ in order to describe it and compute then \eqref{Ereg} by means of the pixels in $\Omega$ and their values of intensity (see Figure \ref{fig:line_pixel}).  It should be noticed that the subdivision curve $\br(t)$ is represented as a polygon with vertices in $\{ \br(i/2^k), i=0,\ldots,2^kM-1 \}$ living on the image. Then, it should be observed that $\br(i/2^k) = (x_i^k,y_i^k)$ is represented on the image by the pixel with coordinates $(\lceil x_i^k \rceil,\lceil y_i^k \rceil)$.

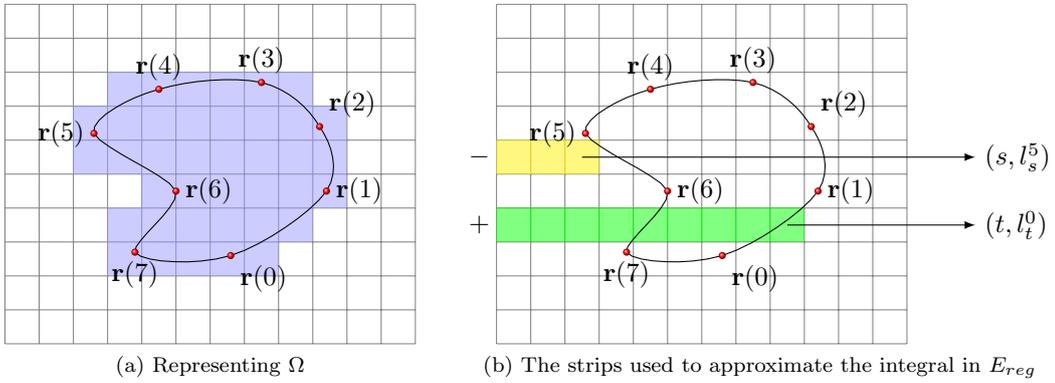
\begin{figure}[h]
  \centering
\subfloat[Representing $\Omega$]{
  \begin{tikzpicture}[scale=.45]
  \draw[help lines] (-2,-2) grid (10,8);

\foreach \i/\j in {0/3,0/4,1/0,1/1,1/3,1/4,1/5,2/0,2/5,3/0,3/5,4/0,4/5,5/0,5/5,6/5,7/2,7/3,7/4}{
  \fill[blue!40,opacity=.5] (\i,\j) -- ++(1,0) -- ++(0,1) -- +(-1,0) -- cycle;
}
\fill[blue!40,opacity=.5] (2,1) -- ++(5,0) -- ++(0,4) -- +(-5,0) -- cycle;

\foreach \x/\y [count=\i] in {4.6/.6, 7.4/2.5, 7.2/4.4, 5.5/5.7, 2.5/5.5, .6/4.2, 3/2.5, 1.8/.7}%
{
\coordinate (P\i) at (\x,\y);
\coordinate (p\i) at (\x,\y);
}
\node[below right] at (p1) {$\br(0)$};
\node[right] at (p2) {$\br(1)$};
\node[above right] at (p3) {$\br(2)$};
\node[above] at (p4) {$\br(3)$};
\node[above] at (p5) {$\br(4)$};
\node[left] at (p6) {$\br(5)$};
\node[right] at (p7) {$\br(6)$};
\node[below] at (p8) {$\br(7)$};
\pgfmathtruncatemacro{\count}{\i}

\ssfourpoints(\i,3)

\draw[solid] (p1) \foreach \j in {2,...,\i}{ -- (p\j)} -- cycle;	

\foreach \j in {1,2,...,\count}{
\shade[ball color = red] (P\j) circle (3pt);
}
\end{tikzpicture}
}
\quad
\subfloat[The strips used to approximate the integral in $E_{reg}$]{
  \begin{tikzpicture}[scale=.45]
  \draw[help lines] (-2,-2) grid (10,8);

\foreach \i in {-2,...,0}
\fill[yellow!100,opacity = .5] (\i,3) -- ++(1,0) -- ++(0,1) -- +(-1,0) -- cycle;
\node at (-2.5, 3.5) {$-$};
\draw[tip] (.5,3.5) -- (12,3.5) node[right] {$(s,l^{5}_s)$};
\foreach \i in {-2,...,6}
\fill[green!100,opacity = .5] (\i,1) -- ++(1,0) -- ++(0,1) -- +(-1,0) -- cycle;
\node at (-2.5, 1.5) {$+$};
\draw[tip] (6.5,1.5) -- (12,1.5) node[right] {$(t,l^{0}_t)$};

\foreach \x/\y [count=\i] in {4.6/.6, 7.4/2.5, 7.2/4.4, 5.5/5.7, 2.5/5.5, .6/4.2, 3/2.5, 1.8/.7}%
{
\coordinate (P\i) at (\x,\y);
\coordinate (p\i) at (\x,\y);
}
\node[below right] at (p1) {$\br(0)$};
\node[right] at (p2) {$\br(1)$};
\node[above right] at (p3) {$\br(2)$};
\node[above] at (p4) {$\br(3)$};
\node[above] at (p5) {$\br(4)$};
\node[left] at (p6) {$\br(5)$};
\node[right] at (p7) {$\br(6)$};
\node[below] at (p8) {$\br(7)$};
\pgfmathtruncatemacro{\count}{\i}

\ssfourpoints(\i,3)

\draw[solid] (p1) \foreach \j in {2,...,\i}{ -- (p\j)} -- cycle;	

\foreach \j in {1,2,...,\count}{
\shade[ball color = red] (P\j) circle (3pt);
}
\end{tikzpicture}
}
  \caption{Pixels considered to compute the {\em region-based energy}. Observe that pixel's coordinates are in the coordinate system of the image {\em (row, column)}.}
  \label{fig:pixels_reg_energy}
\end{figure}

According to our rasterization algorithm of the snake (detailed in {\em Section \ref{sec:implem}}), the integral of the intensity may be approximately computed summing up (with sign) the contribution of each horizontal image strip intersected by $\Omega$, see Figure \ref{fig:pixels_reg_energy}. The value $l_j^i$ is the index of the column of the pixel that results from the intersection of the edge
$[\br\left(i/2^k\right),\br\left((i+1)/2^k\right)]$  with the $j$-th row of the image. Then,

\begin{equation}
\label{eq:int_O_discrete}
I_\Omega = \int\int_\Omega I(x,y) \; \ dx \ dy \approx \sum_{i=0}^{2^kM-1} \text{sign}(x_i^k-x_{i+1}^k) \sum_{j=\left\lceil x_i^k\right\rceil}^{\left\lceil x_{i+1}^k\right\rceil}\sum_{l=1}^{l_j^i} I(l,j).
\end{equation}
The values of $\displaystyle\sum_{l=1}^{l_j^i} I(l,j)$ can be pre-computed in a lookup table to speed up the implementation.

In particular, the approximation of the area of $\Omega$ enclosed by the subdivision curve is,
\begin{equation}
\label{eq:int_areaO_discrete}
\lvert \Omega \rvert = \int\int_\Omega  \; \ dx \ dy
	\approx \sum_{i=0}^{2^kM-1} \text{sign}(x_i^k-x_{i+1}^k) \sum_{j=\left\lceil x_i^k\right\rceil}^{\left\lceil x_{i+1}^k\right\rceil} l_j^i.
\end{equation}

Finally, both approximations \eqref{eq:int_O_discrete} and \eqref{eq:int_areaO_discrete} are substituted in \eqref{Ereg1} to provide the approximation of the region energy.

\subsection{Optimization}
\label{sec:optim}

To obtain the optimal position of the control points of the snake we minimize the total energy given by,
\begin{equation}
E_{\text{snake}}(\bP^0) = \alpha E_{grad}(\bP^0) + (1-\alpha) E_{reg}(\bP^0).
\label{Etotal}
\end{equation}

The optimization problem is solved using the BFGS Quasi-Newton method with a cubic line search procedure. This method requires the gradient of the snake energy with respect to the variables of our problem: the coordinates $(x_j^0,y_j^0)$ of the control points ${\bp}_j^0$, $j=0,\ldots,M-1$.
In this section we give the expressions of the approximations of partial derivatives of each energy with respect to each coordinate $x_j^0$ and $y_j^0$.

\subsubsection{Derivatives of gradient energy}

From \eqref{Egrad2} we obtain (see {\em Appendix \ref{sec:grad_details}} for more details),
\begin{equation}
\frac{\partial E_{grad} }{\partial x_j^0}=\int_{0}^M \left(\left(\frac{\partial^2 I}{\partial x \partial y}\frac{\partial x}{\partial t}-\frac{\partial^2 I}{\partial x^2}\frac{dy}{dt}\right)\frac{\partial x}{\partial x_j^0}+\frac{\partial I}{\partial y}\frac{\partial \left(\frac{dx}{dt}\right) }{\partial x_j^0} \right) \; dt.
\label{derEgradxjk}
\end{equation}

Substituting the integral in \eqref{derEgradxjk} by the average of the integrand evaluated in the parameter values $i/2^k,\;\;i=0,\ldots,2^kM-1$ we obtain (see more details in {\em Appendix \ref{sec:grad_details}}) the following approximation for the partial derivative of gradient energy with respect to $x_j^0$,
\begin{multline} \label{derEgradxjD}
\frac{\partial E_{grad}}{\partial x_j^0} \approx \frac{1}{2^kM}\sum_{i=0}^{2^kM-1} \left(\frac{\partial^2 I}{\partial x \partial y}\left(\br\left(\frac{i}{2^k}\right)\right) t_{ix}^k - \frac{\partial^2 I}{\partial x^2}\left(\br\left(\frac{i}{2^k}\right)\right) t_{iy}^k \right) \varphi\left(\frac{i}{2^k}-j\right) + \\ \frac{1}{2^k M}\sum_{i=0}^{2^kM-1} \frac{\partial I}{\partial y}\left(\br\left(\frac{i}{2^k}\right)\right) \varphi^{'}\left(\frac{i}{2^k}-j\right).
\end{multline}

Proceeding in a similar way, from \eqref{Egrad2} we obtain,
\begin{equation}
\frac{\partial E_{grad} }{\partial y_j^0} = \int_{0}^M \left(\left(\frac{\partial^2 I}{\partial y^2 }\frac{\partial x}{\partial t} - \frac{\partial^2 I}{\partial x \partial y}\frac{dy}{dt}\right) \frac{\partial y}{\partial y_j^k} + \frac{\partial I}{\partial x} \frac{\partial \left(\frac{dy}{dt}\right) }{\partial y_j^k} \right)\; \dt.
\label{derEgradyjk}
\end{equation}

Discretizing the integral with the same procedure, from \eqref{derEgradyjk} we obtain  the following approximation for the partial derivative of gradient energy with respect to $y_j^0$,
\begin{multline} \label{derEgradyjD}
\frac{\partial E_{grad}}{\partial y_j^0} \approx \frac{1}{2^kM} \sum_{i=0}^{2^kM-1} \left(\frac{\partial^2 I}{\partial y^2}\left(\br\left(\frac{i}{2^k}\right)\right) t_{ix}^k - \frac{\partial^2 I}{\partial x \partial y}\left(\br\left(\frac{i}{2^k}\right)\right) t_{iy}^k \right) \varphi\left(\frac{i}{2^k}-j\right) + \\
\frac{1}{2^k M} \sum_{i=0}^{2^kM-1} \frac{\partial I}{\partial x}\left(\br\left(\frac{i}{2^k}\right)\right) \varphi^{'}\left(\frac{i}{2^k}-j\right).
\end{multline}

\subsubsection{Derivatives of region energy}

In order to find the optimal control polygon we have to compute the partial derivatives of $E_{reg}$ with respect to the coordinates $(x_j^0,y_j^0)$ of the control points $\{ \bp_j^0, j=0,\ldots,M-1 \}$.
Since $I_{R}$ is constant, from \eqref{Ereg1} we obtain,
\begin{equation}
\frac{\partial E_{reg}}{\partial x_j^0}=-2D\left(\frac{\partial A}{\partial x_j^0} - \frac{\partial B}{\partial x_j^0} \right)
\label{parcialEreg}
\end{equation}
where
\begin{equation*}
A := \frac{I_{\Omega}}{|\Omega|},  \quad
B := \frac{I_R - I_{\Omega}}{\lvert R \rvert - \lvert \Omega \rvert} \quad \text{and} \quad
D :=  A-B.
\end{equation*}

Proceeding as it is shown in {\em Appendix \ref{sec:grad_details}} it can be proved that \eqref{parcialEreg} is equals to,
\begin{equation}
\frac{\partial E_{reg}}{\partial x_j^0} = -2D \int_{0}^M \left(G - H\,I(\br(t)) \right)\varphi(t-j)y'(t) \dt
\label{parcialEregx}
\end{equation}
where
\begin{equation*}
G := \frac{I_{\Omega}}{\lvert \Omega \rvert^2}+\frac{I_R-I_{\Omega}}{(\lvert R \rvert-\lvert \Omega \rvert)^2} \quad \text{and} \quad
H := \frac{1}{\lvert \Omega \rvert}+\frac{1}{\lvert R \rvert-\lvert \Omega \rvert}.
\end{equation*}

Proceeding in a similar way and deriving in the second equality of \eqref{Green} it is easy to check that
\begin{equation}
\frac{\partial E_{reg}}{\partial y_j^0} = 2D \int_{0}^M \left(G - H\,I(\br(t)) \right) \varphi(t-j)x'(t) \dt.
\label{parcialEregy}
\end{equation}

In practice, we approximate \eqref{parcialEregx} and \eqref{parcialEregy} by
\begin{align*}
\frac{\partial E_{reg}}{\partial x_j^0} &\approx -\frac{\widetilde{D}}{2^{k-1}M} \sum_{i=0}^{2^kM-1}\left[\widetilde{G} - \widetilde{H}I\left(\br\left(\frac{i}{2^k}\right)\right)\right] \varphi\left(\frac{i}{2^k}-j\right)y'\left(\frac{i}{2^k}\right)   \\
\frac{\partial E_{reg}}{\partial y_j^0} &\approx \frac{\widetilde{D}}{2^{k-1}M} \sum_{i=0}^{2^kM-1}\left[\widetilde{G}-\widetilde{H} I\left(\br\left(\frac{i}{2^k}\right)\right)\right] \varphi\left(\frac{i}{2^k}-j\right)x'\left(\frac{i}{2^k}\right)
\end{align*}
where $\widetilde{D},\widetilde{G}$ and $\widetilde{H}$ denote the approximations of $D,G$ and $H$ respectively obtained from
the approximated values of $I_{R},I_{\Omega},\lvert \Omega \rvert$ and $\lvert R \rvert$ in \eqref{eq:int_O_discrete} and \eqref{eq:int_areaO_discrete}.

\section{Implementation}
\label{sec:implem}

In this section we give some details about the computation of the energies previously introduced. Moreover, we describe the main features of the application SubdivisionSnake, which is able to compute the subdivision snakes produced by  cubic B-spline and 4-point subdivision curves.

\subsection{Details about the energies}

For the implementation of energies it's necessary to define how to compute the gradient of an image in a point, the area enclosed by a curve, and others details. In the following we discuss these themes.

\subsubsection{Gradient energy}

Since the image $I$ is only defined in points with integer coordinates, the evaluation of $I$ and its partial derivatives
in a point $(x,y)\in \RR^2$ is approximated using bilinear interpolation. In particular, $\nabla I$ in \ref{Egrad3} is approximated as,

\begin{multline}
\nabla I(x,y) := \nabla I(\lfloor x \rfloor, \lfloor y \rfloor) \,(1-\{x\})(1-\{y\}) +
\nabla I(\lfloor x+1 \rfloor, \lfloor y \rfloor) \,\{x\}(1-\{y\}) + \\
\nabla I(\lfloor x \rfloor, \lfloor y+1 \rfloor) \,(1-\{x\})\{y\} +
\nabla I(\lfloor x+1 \rfloor, \lfloor y+1 \rfloor)\{x\}\{y\},
\end{multline}
where $\{x\} = x - \lfloor x \rfloor$ denotes the fractional part of $x$.

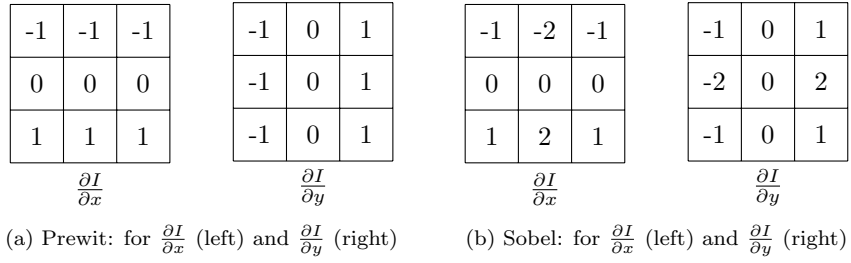
\begin{figure}[h]
\centering
\subfloat[Prewit: for $\frac{\partial I}{\partial x}$ (left) and $\frac{\partial I}{\partial y}$ (right)]{
\begin{tikzpicture}[scale=.7]
\draw (0,0) grid (3,3);
\foreach \i in {1,2,3}{
\node at (-.5 +\i,.5) {1};
\node at (-.5 +\i,1.5) {0};
\node at (-.5 +\i,2.5) {-1};
}
\node at (1.5,-.5) {$\frac{\partial I}{\partial x}$};
\end{tikzpicture}
\qquad
\begin{tikzpicture}[scale=.7]
\draw (0,0) grid (3,3);
\foreach \i in {1,2,3}{
\node at (.5,-.5 +\i) {-1};
\node at (1.5,-.5 +\i) {0};
\node at (2.5,-.5 +\i) {1};
}

\node at (1.5,-.5) {$\frac{\partial I}{\partial y}$};
\end{tikzpicture}
}
\qquad
\subfloat[Sobel: for $\frac{\partial I}{\partial x}$ (left) and $\frac{\partial I}{\partial y}$ (right)]{
\begin{tikzpicture}[scale=.7]
\draw (0,0) grid (3,3);
\foreach \i in {1,3}{
\node at (-.5 +\i,.5) {1};
\node at (-.5 +\i,2.5) {-1};
}
\foreach \i in {1,2,3} \node at (-.5 +\i,1.5) {0};

\node at (1.5,.5) {2};
\node at (1.5,2.5) {-2};
\node at (1.5,-.5) {$\frac{\partial I}{\partial x}$};
\end{tikzpicture}
\qquad
\begin{tikzpicture}[scale=.7]
\draw (0,0) grid (3,3);
\foreach \i in {1,3}{
\node at (.5,-.5 +\i) {-1};
\node at (2.5,-.5 +\i) {1};
}
\foreach \i in {1,2,3} \node at (1.5,-.5 +\i) {0};

\node at (.5,1.5) {-2};
\node at (2.5,1.5) {2};
\node at (1.5,-.5) {$\frac{\partial I}{\partial y}$};
\end{tikzpicture}
}
\caption{Some filters known to compute approximations of the gradient in a pixel.}
\label{fig:filters}
\end{figure}

\begin{figure}[h]
\centering
\subfloat[Image $I$]{
\includegraphics[width=.19\textwidth]{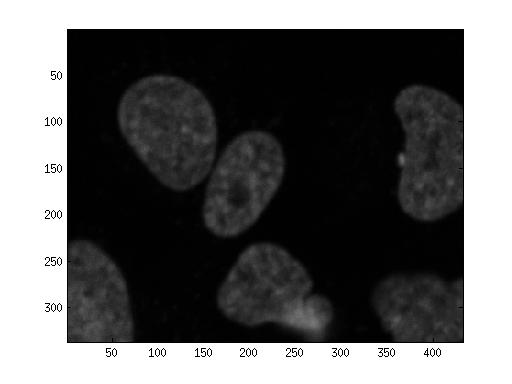}
}
\subfloat[Filter Sobel: $\frac{\partial I}{\partial x}$]{
\includegraphics[width=.19\textwidth]{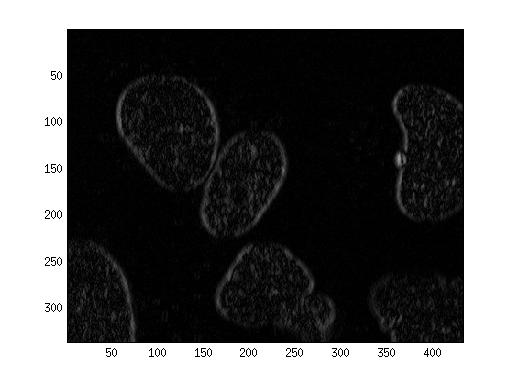}
}
\subfloat[Filter Sobel: $\frac{\partial I}{\partial y}$]{
\includegraphics[width=.19\textwidth]{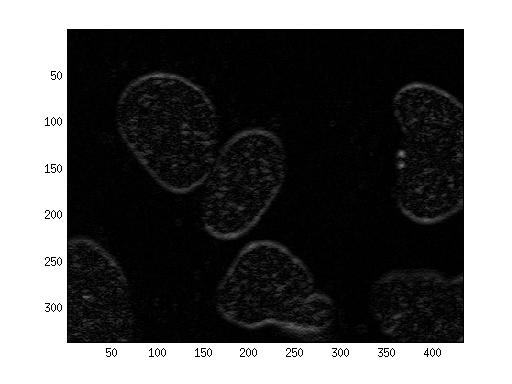}
}
\subfloat[Filter $7\times7$: $\frac{\partial I}{\partial x}$]{
\includegraphics[width=.19\textwidth]{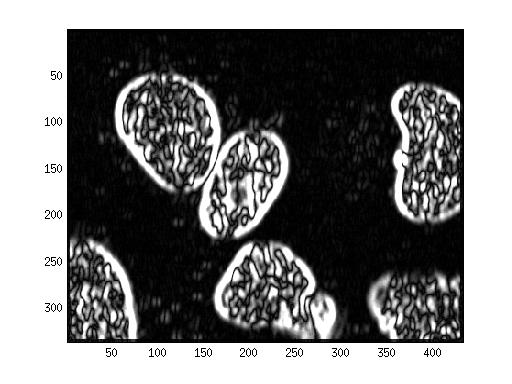}
}
\subfloat[Filter $7\times7$: $\frac{\partial I}{\partial y}$]{
\includegraphics[width=.19\textwidth]{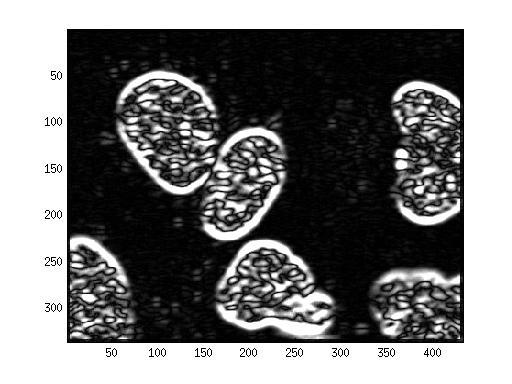}
}
\caption{Partial derivatives computed using filters of different sizes. a) Original images, b) and c) Sobel's filter, d) and e) the proposed $7 \times 7$ filter.}
\label{fig:image_filter}
\end{figure}

The gradient of the image in a pixel can be approximated using different filters such as Prewitt and Sobel \cite{Gonzo} (see Figure \ref{fig:filters}). Since we evaluate the gradient in points that belong to the snake, it is convenient to extend the width of the filter in order to increase the region of attraction of gradient energy (see Figure \ref{fig:image_filter}). Consequently, we use a generalization of the Prewitt filter of ${(2q+1)}\times{(2q+1)}$ pixels, to compute the gradient in those pixels with distance greater or equal to  $q>0$ (see Figure \ref{fig:prop_GradFilter}) to the boundary of the image. For the rest of the pixels we use Sobel filter to approximate the gradient. The constant value $q$ depends on the image dimensions.

\begin{figure}[h]
\centering
\begin{tikzpicture}[scale=.6]
\draw (0,2) grid +(2,3);
\foreach \i in {1,2}{
\node at (-.5 +\i,2.5) {1};
\node at (-.5 +\i,3.5) {0};
\node at (-.5 +\i,4.5) {-1};
}
\foreach \i in {0,1,2}
\node at (2.5,.5 +3*\i) {\ldots};

\draw (3,2) grid +(2,3);
\foreach \i in {1,2}{
\node at (2.5 +\i,2.5) {1};
\node at (2.5 +\i,3.5) {0};
\node at (2.5 +\i,4.5) {-1};
}

\foreach \x/\y in {1/1.5,4/1.5,1/5.5,4/5.5}
\node at (\x,\y) {\vdots};

\foreach \x/\y in {0/0,3/0}{
\draw (\x,\y) grid +(2,1);
\node at (\x +.5,.5) {1};
\node at (\x +1.5,.5) {1};
}
\foreach \x/\y in {0/6,3/6}{
\draw (\x,\y) grid +(2,1);
\node at (\x +.5,6.5) {-1};
\node at (\x +1.5,6.5) {-1};
}

\foreach \i in {1,5}
\node at (2.5,.5 + \i) {$\ddots$};

\draw[decorate,decoration={brace,raise=5pt},thick] (0,0) -- node[left=7pt] {$q$} (0,3);
\draw[decorate,decoration={brace,raise=5pt},thick] (0,4) -- node[left=7pt] {$q$} (0,7);
\draw[decorate,decoration={brace,raise=5pt},thick] (0,7) -- node[above=7pt] {$2q+1$} (5,7);
\end{tikzpicture}
\qquad \qquad
\begin{tikzpicture}[scale=.6]
\draw (2,0) grid +(3,2);
\foreach \i in {1,2}{
\node at (2.5,-.5 +\i) {-1};
\node at (3.5,-.5 +\i) {0};
\node at (4.5,-.5 +\i) {1};
}
\foreach \i in {0,1,2}
\node at (.5 +3*\i,2.5) {\vdots};

\draw (2,3) grid +(3,2);
\foreach \i in {1,2}{
\node at (2.5,2.5 +\i) {-1};
\node at (3.5,2.5 +\i) {0};
\node at (4.5,2.5 +\i) {1};
}

\foreach \x/\y in {1.5/1,1.5/4,5.5/1,5.5/4}
\node at (\x,\y) {\ldots};

\foreach \x/\y in {0/0,0/3}{
\draw (\x,\y) grid +(1,2);
\node at (.5,\y +.5) {-1};
\node at (.5,\y +1.5) {-1};
}
\foreach \x/\y in {6/0,6/3}{
\draw (\x,\y) grid +(1,2);
\node at (6.5,\y +.5) {1};
\node at (6.5,\y +1.5) {1};
}

\foreach \i in {1,5}
\node at (.5 + \i,2.5) {$\ddots$};

\draw[decorate,decoration={brace,raise=5pt},thick] (0,5) -- node[above=7pt] {$q$} (3,5);
\draw[decorate,decoration={brace,raise=5pt},thick] (4,5) -- node[above=7pt] {$q$} (7,5);
\draw[decorate,decoration={brace,raise=5pt},thick] (0,0) -- node[left=7pt] {$2q+1$} (0,5);
\end{tikzpicture}
\caption{Proposed filter to compute approximations of the gradient of an image:  for $\frac{\partial I}{\partial x}$ (left) and $\frac{\partial I}{\partial y}$ (right).}
\label{fig:prop_GradFilter}
\end{figure}

The gradient of the image in each pixel is precomputed and stored in a lookup table, so that the evaluations in \eqref{Egrad3}, \eqref{derEgradxjD} and \eqref{derEgradyjD} use the stored values.

\subsubsection{Region energy}

The first step to compute the integrals \eqref{Ereg} defining the region energy is to obtain a sequence of pixels that approximates the snake, that is represented by the polygon with vertices $\{\br(i/2^k) = (x_i^k,y_i^k), i=0,\ldots,2^kM-1 \}$. The problem is reduced to the rasterization of each edge of that polygon. Rasterization algorithms provide the pixels that are intersected by a straight line (see Figure \ref{subfig:rasterization}). Since these are more pixels than the ones needed to describe the region $\Omega$ enclosed by a closed polygon, we select for each horizontal line only one pixel per edge of the polygon \footnote{We choose the horizontal direction without loss of generality, the same result is obtained if the vertical direction is chosen.}. To obtain these pixels, called {\em boundary pixels}, we determine for the horizontal line $j$ the pixels $(j,l^i_j)$  that are simultaneously on the line and on the edge $[\br\left(i/2^k\right),\br\left((i+1)/2^k\right)]$. If the result of the previous operation is more than one pixel, then we select the outer pixel with respect to the region enclosed by the subdivision curve (see Figure \ref{subfig:boundary}). Consequently, we proceed as follows.

We classify the edge $[\br\left(i/2^k\right),\br\left((i+1)/2^k\right)]$ as {\it downhill, horizontal} or {\it uphill} if the sign of $x_{i}^k - x_{i+1}^k$ is negative, zero or positive, respectively \footnote{Remember that we are using the system of coordinates defined by $(row,column)$.}.
To compute approximately the integrals in \eqref{Ereg} it is necessary to chose, for a given edge $[\br\left(i/2^k\right),\br\left((i+1)/2^k\right)]$, one pixel with coordinates $(j,l_j^i)$ for each image row $j$, with $\min\{\lceil x_i^k \rceil,\lceil x_{i+1}^k \rceil\} \leq j \leq  \max\{\lceil x_i^k \rceil,\lceil x_{i+1}^k \rceil\}$. The value of $l_j^i$ depends on the previous edge classification as follows. Let $r_i(x)$ be the equation of the line passing through the pixels $(\lceil x_i^k \rceil,\lceil y_i^k \rceil)$ and $(\lceil x_{i+1}^k \rceil,\lceil y_{i+1}^k \rceil)$, then,
\[
r_i(x) = \lceil y_i^k \rceil + \frac{\lceil y_{i+1}^k \rceil - \lceil y_i^k \rceil}{\lceil x_{i+1}^k \rceil - \lceil x_i^k \rceil}(x- \lceil x_i^k \rceil).
\]

If $[\br\left(i/2^k\right),\br\left((i+1)/2^k\right)]$ is a downhill edge (see Figure \ref{subfig:boundary}, edge $[\br\left(4\right),\br\left(5\right)]$) then,
\begin{equation}
  \label{eq:line_pixel_lowerD}
  l_j^i =
    \min \left\{ \left\lceil r_i(j) \right\rceil , \left\lceil r_i(j+1) \right\rceil \right\}, \qquad j \in \left[ \left\lceil x_i^k \right\rceil,  \left\lceil x_{i+1}^k \right\rceil\right]
\end{equation}

\begin{figure}[h]
\centering
\subfloat[Rasterization of straight lines]{\label{subfig:rasterization}
\begin{tikzpicture}[scale=.7]
\draw[help lines] [step=0.5cm] (-1,-1.5) grid (8.5,7.5);


\foreach \x/\y [count=\i] in {4.6/.6, 7.1/2.7, 5.5/5.7, 2.5/6.2, .65/4.25, 3.1/2.6, 1.8/.7}%
{
\coordinate (P\i) at (\x,\y);
\coordinate (p\i) at (\x,\y);
}
\foreach \i/\j in {2.5/2.5,3/2.5, 1.5/3,2/3, 1/3.5,1.5/3.5, .5/4}{
  \fill[blue!40,opacity=.8] (\i,\j) -- ++(.5,0) -- ++(0,.5) -- +(-.5,0) -- cycle;
}
\foreach \i/\j in {6.5/2.5,7/2.5, 6.5/3, 6/3.5,6.5/3.5, 6/4, 5.5/4.5,6/4.5, 5.5/5,5.5/5.5}{
  \fill[blue!40,opacity=.8] (\i,\j) -- ++(.5,0) -- ++(0,.5) -- +(-.5,0) -- cycle;
}

\node[below right] at (p1) {$\br(0)$};
\node[right] at (p2) {$\br(1)$};
\node[above right] at (p3) {$\br(2)$};
\node[above] at (p4) {$\br(3)$};
\node[left] at (p5) {$\br(4)$};
\node[right] at (p6) {$\br(5)$};
\node[left] at (p7) {$\br(6)$};
\pgfmathtruncatemacro{\count}{\i}

\ssfourpoints(\i,3)

\draw[dotted,very thick,opacity=.5] (p1) \foreach \j in {2,...,\i}{ -- (p\j)} -- cycle;	
\draw[very thick] (p9) -- (p17);
\draw[very thick] (p33) -- (p41);

\foreach \j in {1,2,...,\count}{
\shade[ball color = red] (P\j) circle (3pt);
}

\end{tikzpicture}
}
\quad
\subfloat[Boundary pixels to describe the boundary]{\label{subfig:boundary}
\begin{tikzpicture}[scale=.7]
\draw[help lines] [step=0.5cm] (-1,-1.5) grid (8.5,7.5);


\foreach \x/\y [count=\i] in {4.6/.6, 7.1/2.7, 5.5/5.7, 2.5/6.2, .65/4.25, 3.1/2.6, 1.8/.7}%
{
\coordinate (P\i) at (\x,\y);
\coordinate (p\i) at (\x,\y);
}
\foreach \i/\j in {2.5/2.5,1.5/3,1/3.5,.5/4}{
  \fill[blue!40,opacity=.8] (\i,\j) -- ++(.5,0) -- ++(0,.5) -- +(-.5,0) -- cycle;
}
\foreach \i/\j in {7/2.5,6.5/3,6.5/3.5,6/4,6/4.5,5.5/5,5.5/5.5}{
  \fill[blue!40,opacity=.8] (\i,\j) -- ++(.5,0) -- ++(0,.5) -- +(-.5,0) -- cycle;
}

\node[below right] at (p1) {$\br(0)$};
\node[right] at (p2) {$\br(1)$};
\node[above right] at (p3) {$\br(2)$};
\node[above] at (p4) {$\br(3)$};
\node[left] at (p5) {$\br(4)$};
\node[right] at (p6) {$\br(5)$};
\node[left] at (p7) {$\br(6)$};
\pgfmathtruncatemacro{\count}{\i}

\ssfourpoints(\i,3)

\draw[dotted,very thick,opacity=.5] (p1) \foreach \j in {2,...,\i}{ -- (p\j)} -- cycle;	
\draw[very thick] (p9) -- (p17);
\draw[very thick] (p33) -- (p41);

\foreach \j in {1,2,...,\count}{
\shade[ball color = red] (P\j) circle (3pt);
}

\end{tikzpicture}
}
\caption{Pixel discretization of a straight line for a left edge and a right edge describing the boundary of a region.}
\label{fig:line_pixel}
\end{figure}

If $[\br\left(i/2^k\right),\br\left((i+1)/2^k\right)]$ is a uphill edge (see Figure \ref{subfig:boundary}, edge $[\br\left(1\right),\br\left(2\right)]$) then,
\begin{equation}
  \label{eq:line_pixel_lowerU}
  l_j^i =
    \max \left\{ \left\lceil r_i(j) \right\rceil , \left\lceil r_i(j+1) \right\rceil \right\}, \qquad j \in \left[ \left\lceil x_{i+1}^k \right\rceil, \left\lceil x_{i}^k \right\rceil\right]
\end{equation}

Finally, if $[\br\left(i/2^k\right),\br\left((i+1)/2^k\right)]$ is a horizontal edge, then there is no need to define the value of $l_j^i$ as $\text{sign}(x_i^k-x_{i+1}^k) = 0$ in \eqref{eq:int_O_discrete} and \eqref{eq:int_areaO_discrete}. In this case, the description of the boundary makes use of the neighboring edges.

\begin{figure}[h]
\centering
\begin{tikzpicture}[scale=.45]
  \draw[help lines] (1,1) grid (13,12);

\foreach \i/\j in {6/10,9/10, 2/9,3/9,6/9,11/9, 2/8,4/8,5/8,11/8, 2/7,5/7,11/7, 2/6,10/6, 2/5,10/5, 3/4,10/4, 3/3,10/3, 5/2,10/2}{
  \fill[blue!40,opacity=.5] (\i,\j) -- ++(1,0) -- ++(0,1) -- +(-1,0) -- cycle;
}
\foreach \i/\j in {7/10,8/10, 7/9,8/9,9/9,10/9, 3/8,6/8,7/8,8/8,9/8,10/8, 3/7,4/7,6/7,7/7,8/7,9/7,10/7, 3/6,4/6,5/6,6/6,7/6,8/6,9/6, 3/5,4/5,5/5,6/5,7/5,8/5,9/5, 4/4,5/4,6/4,7/4,8/4,9/4, 4/3,5/3,6/3,7/3,8/3,9/3, 6/2,7/2,8/2,9/2}{
  \fill[gray!40,opacity=.5] (\i,\j) -- ++(1,0) -- ++(0,1) -- +(-1,0) -- cycle;
}

\foreach \x/\y [count=\i] in {2.6/9.8, 2.5/6.7, 4.6/3.2, 10.4/2.7, 10.5/5.5, 11.4/8.85, 7.2/10.4, 5.4/7.5}%
{
\coordinate (P\i) at (\x,\y);
\coordinate (p\i) at (\x,\y);
}
\node[above] at (p1) {$\br(0)$};
\node[left] at (p2) {$\br(1)$};
\node[below left] at (p3) {$\br(2)$};
\node[below right] at (p4) {$\br(3)$};
\node[right] at (p5) {$\br(4)$};
\node[right] at (p6) {$\br(5)$};
\node[above] at (p7) {$\br(6)$};
\node[below] at (p8) {$\br(7)$};
\pgfmathtruncatemacro{\count}{\i}

\ssfourpoints(\i,3)

\draw[solid] (p1) \foreach \j in {2,...,\i}{ -- (p\j)} -- cycle;	

\foreach \j in {1,2,...,\count}{
\shade[ball color = red] (P\j) circle (3pt);
}
\end{tikzpicture}
\caption{Description of the boundary of a region with pairs of {\em boundary pixels}.}
\label{fig:convexity_bound}
\end{figure}
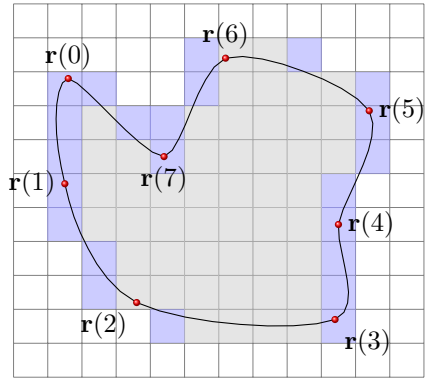

In order to describe the boundary of the region enclosed by the subdivision curve, we store pairs of {\em boundary pixels} with respect to each horizontal. The amount of pairs on each horizontal line depends on the convexity of the curve (see Figure \ref{fig:convexity_bound}). It should be noticed that a boundary pixel may be simultaneously the right pixel of a pair and left pixel of the next pair in the same horizontal line, see for example the pixel corresponding to $\br(7)$ in Figure \ref{fig:convexity_bound}.

\subsection{{\bf SubdivisionSnake} Application}

The segmentation technique proposed in this paper has been implemented in C\# with \texttt{.NET} platform (version 4.5). The optimization step, based on Limited memory BFGS method \cite{Byrd95}, is done using the library  \texttt{Optimization.dll} of \texttt{Accord.NET} \cite{ANET15}. The resulting application is called {\bf SubdivisionSnake} and it is independent of any imaging hardware.

\begin{figure}[!ht]\centering
\subfloat[Initialization]{
\includegraphics[width=.3\textwidth]{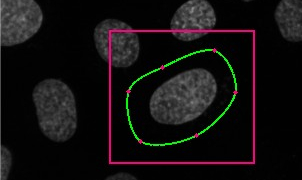}
}
\subfloat[Segmentation]{
\includegraphics[width=.3\textwidth]{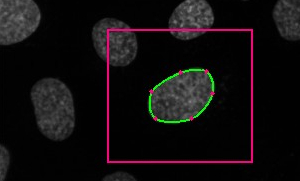}
}
\subfloat[Interaction with the control points]{
\includegraphics[width=.3\textwidth]{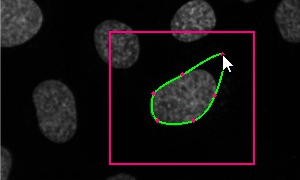}
}
\caption{Interaction of the user with SubdivisionSnake application.}
\label{fig:snake1}
\end{figure}

Currently, {\bf SubdivisionSnake} application is able to compute two type of subdivision snakes: cubic B-spline \eqref{BS_evenrule}-\eqref{BS_oddrule} and four points subdivision curves \eqref{4pt_evenrule}-\eqref{4pt_oddrule}. Other subdivision curves could be easily added to the application if a procedure for evaluating the curve and its partial derivatives at dyadic parametric values is included. The interaction with the user is very simple and only requires an initial polygon approximating the boundary of the object to be segmented and a bounding box containing the object to be segmented and the initial polygon. As illustrated in Figure \ref{fig:snake1}, the position of any control point can be intuitively manipulated on the image with simple mouse actions. The snake is updated in real-time since control points have local influence and therefore only a small region of the snake has to be recomputed. The resulting tool is a semi-automatic and intuitive segmentation algorithm based on the position of the control points and consisting of
three fundamental steps: initialization, optimization and correction.

{\bf SubdivisionSnake} application has a main window to load the image and set the free parameters.
This window contains the following options:

\begin{itemize}
\item Load and save image: this option load the target image and save the segmented image. The application accepts \texttt{.jpg}, \texttt{.gif} , \texttt{.bmp}, \texttt{.png} and \texttt{.tif} images. Color images are transformed in gray images at the beginning of the processing.

\item Target color object: using this option the user says if the object to be segmented is darker than the background (default option) or the contrary.

\item Control points: the application offers several options related to the control points of the snake.

Clicking on the image, the initial position of control points can be defined. It is also possible to delete a control point or to change its position (dragging the mouse). Finally, control points can be saved or loaded to be reused.

\item $\alpha$: this options allows the selection of a value in $[0,1]$ for the parameter $\alpha$.

By default, we use $\alpha=0.1$ at the beginning of the optimization and $\alpha=0.9$ in the last steps of the optimization since gradient based energies have a narrow zone of attraction in comparison with region based
energies. This means that in the default option the region energy controls initially the movement of the snake inducing his fast displacement. When the position of the snake stabilizes the value of $\alpha$ changes  automatically to $\alpha=0.9$ and then the gradient energy pushes the snake to the boundary of the object.

Otherwise, the selected $\alpha$ does not change during all the optimization process.

\item Visualize: this option is used to visualize the snake and its control polygon.
\end{itemize}

\section{Results}
\label{sec:results}

To illustrate the performance of {\bf SubdivisionSnake} we experiment  with two group of images. The first group is composed by synthetic images. These images are created in such a way that the boundary of the object to be segmented is clear e intuitive. Some synthetic images used in this paper
were produced filling the interior of closed subdivision curves with a color that makes a good contrast with the background. Other synthetic images were obtained using \texttt{.seg} files of Berkeley data base \cite{Berkeley}. From these files it is possible to know which pixels belong to each object to be segmented in the image. The syntectic image is obtained assigning a specific color to these pixels and a contrasting color to the rest of pixels. In general, synthetic images are simpler than the real images included in the second group.
In our experiment we also use real images from Berkeley data base and from other sources. All the examples shown in this section are the direct result of the
optimization process, without any interactive correction. Color images are transformed to gray level images using the standard perceptual weightings for the three-color components \cite{Gonzo}.

\subsection{Quantitative evaluation of results}

When we work with synthetic images the ground-truth region, composed by pixels belonging to the object $\Gamma$, is known.
In some real images of  Berkeley database, the ground-truth is also given. In all these cases it is possible to validate quantitatively the quality of the results using the Jaccard distance $J$ between $\Gamma$ and the region $\Omega$ enclosed by the snake, given by
\[
J = 1-\dfrac{\lvert\Omega \cap \Gamma\rvert}{\lvert\Omega \cup \Gamma\rvert}
\]
where $\vert G \rvert$ denotes the area of region $G$.  Observe that $0 \leq J \leq 1 $ and a value of $J$ close to 0 indicates
a good segmentation.

\begin{figure}[!ht]\centering
\subfloat[Initialization ]{
\includegraphics[width=.24\textwidth]{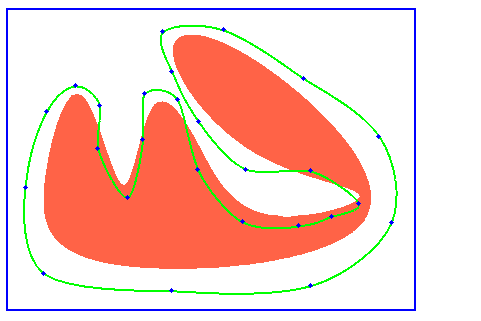}
}
\subfloat[Segmentation ]{
\includegraphics[width=.24\textwidth]{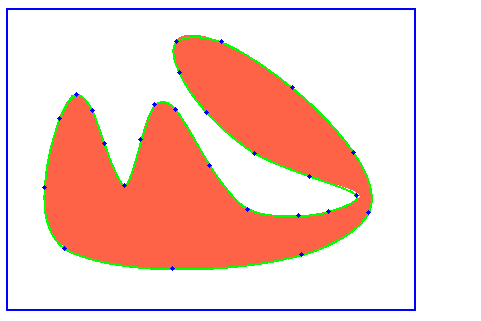}
}
\subfloat[Initialization]{
\includegraphics[width=.24\textwidth]{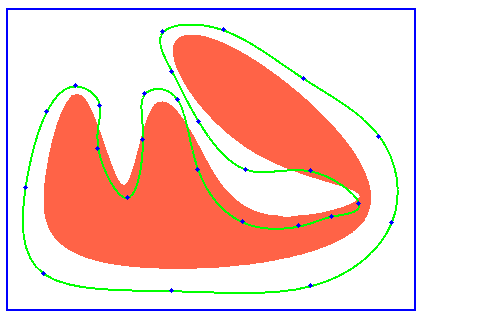}
}
\subfloat[Segmentation ]{
\includegraphics[width=.24\textwidth]{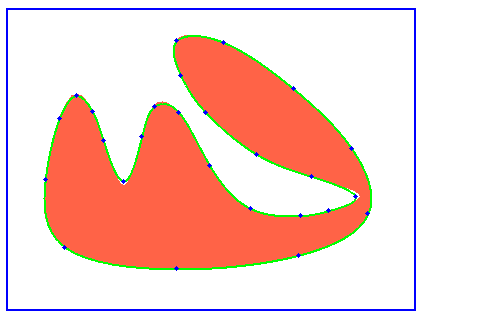}
}
\caption{Syntethic image corresponding to the results reported in Table 1. a) and b) 4 points snake, c) and d) cubic B-spline snake.}
\label{fig:imagesJaccard_artificial}
\end{figure}

\begin{figure}[!ht]\centering
\subfloat[Initialization ]{
\includegraphics[width=.24\textwidth]{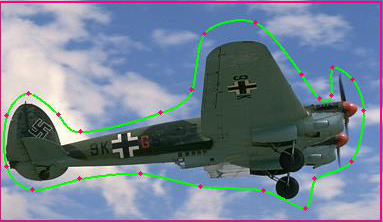}
}
\subfloat[Segmentation ]{
\includegraphics[width=.24\textwidth]{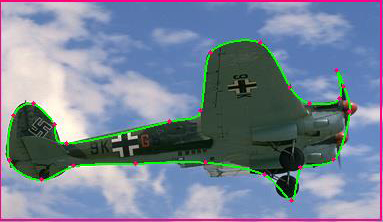}
}
\subfloat[Initialization ]{
\includegraphics[width=.24\textwidth]{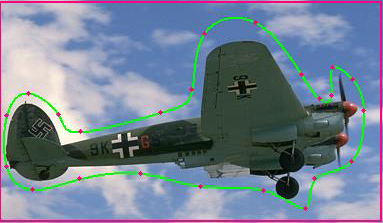}
}
\subfloat[Segmentation ]{
\includegraphics[width=.24\textwidth]{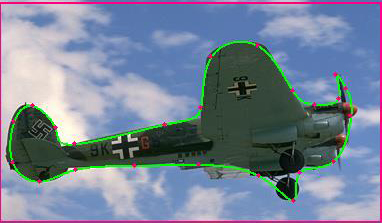}
}
\caption{Airplane image corresponding to the results reported in Table 1. a) and b) 4 points snake, c) and d) cubic B-spline snake.}
\label{fig:imagesJaccard_avion}
\end{figure}

\begin{figure}[!ht]\centering
\subfloat[Initialization ]{
\includegraphics[width=.24\textwidth]{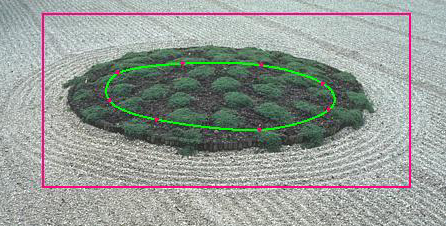}
}
\subfloat[Segmentation ]{
\includegraphics[width=.24\textwidth]{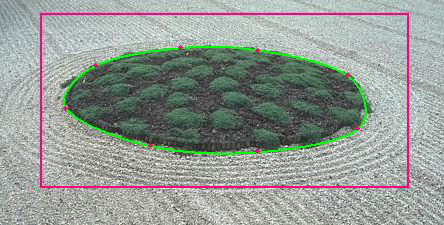}
}
\subfloat[Initialization ]{
\includegraphics[width=.24\textwidth]{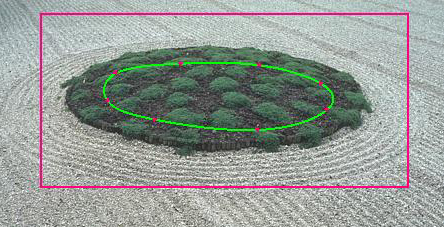}
}
\subfloat[Segmentation ]{
\includegraphics[width=.24\textwidth]{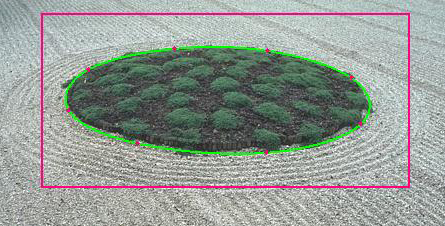}
}
\caption{Japanese garden image corresponding to the results reported in Table 1. a) and b) 4 points snake, c) and d) cubic B-spline snake.}
\label{fig:imagesJaccard_jardin}
\end{figure}

Table \ref{tabla:Jaccard} shows the segmentation results obtained for images in Figures \ref{fig:imagesJaccard_artificial},
\ref{fig:imagesJaccard_avion} and \ref{fig:imagesJaccard_jardin}, with the default selection of the parameter $\alpha$. We use
the same sequence of control points $\mathbf{P}^0$ to initialize the snake based on cubic B-splines and on 4-point subdivision
scheme. In the case of the B-splines, we compute initially the sequence of points $\widetilde{\bP}^0$ \eqref{eq:pseudo_int} such
that the corresponding cubic B-spline interpolates the points $\mathbf{P}^0$. Columns 2 and 3 of Table \ref{tabla:Jaccard} contain
the Jaccard distance between $\Gamma$ and the region $\Omega$ enclosed by the cubic B-spline snake in the initialization step and
after convergence, respectively. Similarly, columns 4 and 5 contain the Jaccard distance for the 4 points subdivision snake.
We observe that despite of the different nature of the images, the Jaccard distance in the optimum is very small for both
subdivision snakes, in correspondence with a good segmentation of the target objet.

\medskip
\begin{table}
\centering
\begin{tabular}{|c|c|c|c|c|c|}  \hline
Image & \# control points & \multicolumn{2}{|c|}{cubic B-spline} & \multicolumn{2}{|c|}{4-point} \\ \cline{3-6}
      &                   &  Initialization & Segmentation     & Initialization & Segmentation \\ \hline
Synthetic image & 26 & 0.4019 & 0.0168 & 0.3935 & 0.0196 \\ \hline
Airplane & 23  & 0.3137 & 0.0837 & 0.3099 & 0.0772 \\ \hline
Japanese garden & 8 & 0.5091 & 0.0471 & 0.5134 & 0.0463 \\ \hline 
\end{tabular}
\caption{Jaccard distance for images in Figures \ref{fig:imagesJaccard_artificial},
\ref{fig:imagesJaccard_avion} and \ref{fig:imagesJaccard_jardin}.}
\label{tabla:Jaccard}
\end{table}

\subsection{Influence of the number of control points}

\begin{figure}[!ht]\centering
\subfloat[Initialization]{
\includegraphics[width=.24\textwidth]{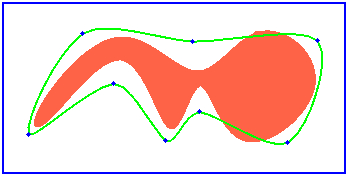}
}
\subfloat[Segmentation]{
\includegraphics[width=.24\textwidth]{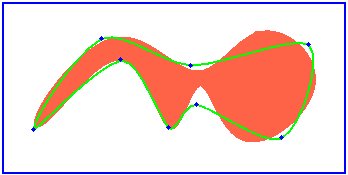}
}
\subfloat[Initialization]{
\includegraphics[width=.24\textwidth]{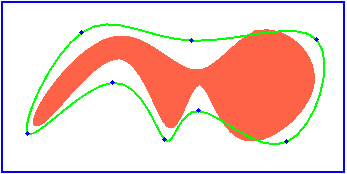}
}
\subfloat[Segmentation]{
\includegraphics[width=.24\textwidth]{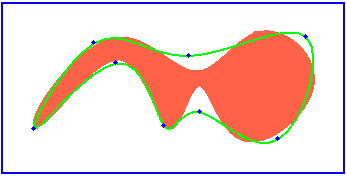}
}
\caption{Images corresponding to the results in Table 2 with 8 control points. a) and b) 4 points snake, c) and d) cubic B-spline snake. }
\label{fig:imagesCP1}
\end{figure}

\begin{figure}[!ht]\centering
\subfloat[Initialization]{
\includegraphics[width=.24\textwidth]{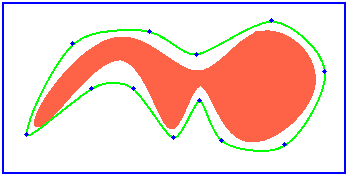}
}
\subfloat[Segmentation]{
\includegraphics[width=.24\textwidth]{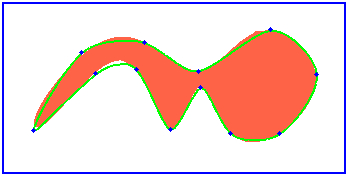}
}
\subfloat[Initialization]{
\includegraphics[width=.24\textwidth]{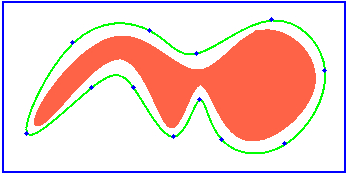}
}
\subfloat[Segmentation]{
\includegraphics[width=.24\textwidth]{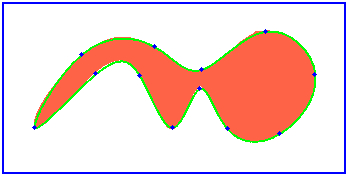}
}
\caption{Images corresponding to the results in Table 2 with 12 control points.  a) and b) 4 points snake, c)and d) cubic B-spline snake. }
\label{fig:imagesCP3}
\end{figure}

The number of control points has a strong influence in the quality of the segmentation.
In general, increasing the number of control points improves the quality of segmentation, but makes higher
the computational cost.  Table \ref{tabla:Jaccard_nP} shows that if one selects initial polygons with different
number of control points but approximately the same Jaccard distance, then the better segmentation corresponds to the snake with the highest number of control points. This result is valid for both subdivision snakes.
Figures \ref{fig:imagesCP1} and \ref{fig:imagesCP3} correspond to the results in Table 2 for 8 and 12 control points respectively. The improvement of the segmentation is evident when we compare Figure \ref{fig:imagesCP1} b) with Figure \ref{fig:imagesCP3} b) and Figure \ref{fig:imagesCP1} d) with Figure \ref{fig:imagesCP3} d).

\begin{table}
\centering
\begin{tabular}{|c|c|c|c|c|}  \hline
\# control points & \multicolumn{2}{|c|}{cubic B-spline} & \multicolumn{2}{|c|}{4-point} \\ \cline{2-5}
      &            Initialization & Segmentation     & Initialization & Segmentation \\ \hline
8   & 0.3630 & 0.1653  & 0.3581  & 0.1873  \\ \hline
10  & 0.3682 & 0.0637  & 0.3539  & 0.0788 \\ \hline
12  & 0.3681 & 0.0388  & 0.3411  & 0.0589 \\ \hline
\end{tabular}
\caption{Influence of the number of control points in the quality of the segmentation measured by the Jaccard distance.}
\label{tabla:Jaccard_nP}
\end{table}

\subsection{Some results with real images}

\begin{figure}[!ht]\centering
\subfloat[Initialization]{
\includegraphics[width=.21\textwidth]{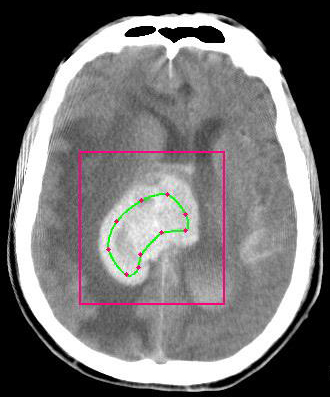}
}
\subfloat[Segmentation]{
\includegraphics[width=.21\textwidth]{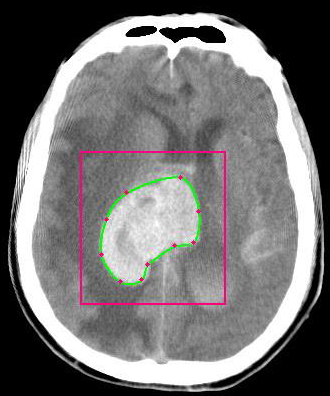}
}
\subfloat[Initialization]{
\includegraphics[width=.21\textwidth]{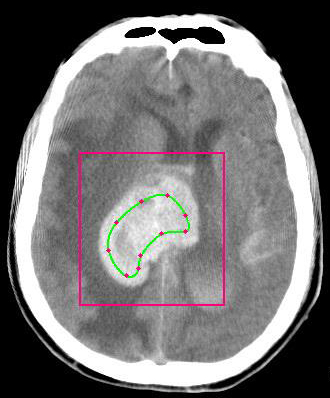}
}
\subfloat[Segmentation]{
\includegraphics[width=.21\textwidth]{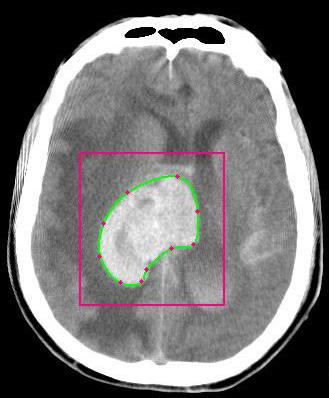}
}
\caption{Segmentation of a cerebral hemorrhage.a) and b) 4 points snake, c) and d) cubic B-spline snake.}
\label{fig:imagesHemorrage}
\end{figure}

\begin{figure}[!ht]\centering
\subfloat[Initialization]{
\includegraphics[width=.24\textwidth]{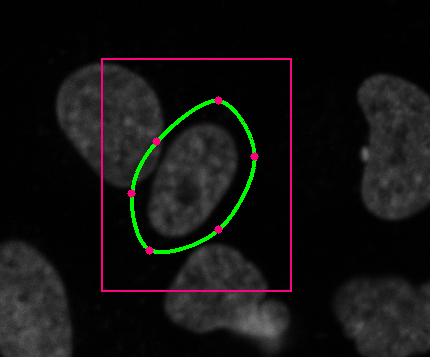}
}
\subfloat[Segmentation]{
\includegraphics[width=.24\textwidth]{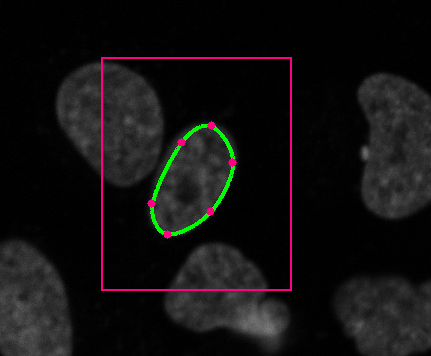}
}
\subfloat[Initialization]{
\includegraphics[width=.24\textwidth]{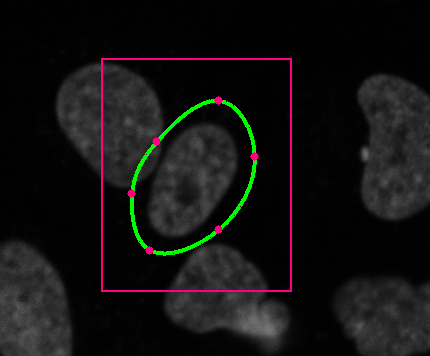}
}
\subfloat[Segmentation]{
\includegraphics[width=.24\textwidth]{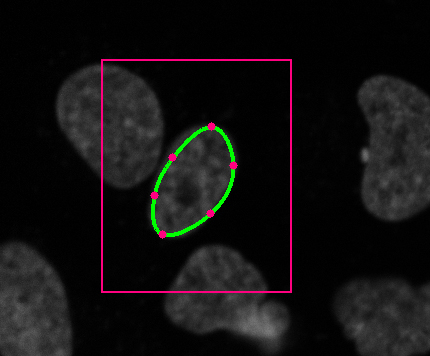}
}
\caption{Segmentation of a cell in a sample. a) and b) 4 points snake, c) and d) cubic B-spline snake.}
\label{fig:imagesCells}
\end{figure}

\begin{figure}[!ht]\centering
\subfloat[Initialization]{
\includegraphics[width=.24\textwidth]{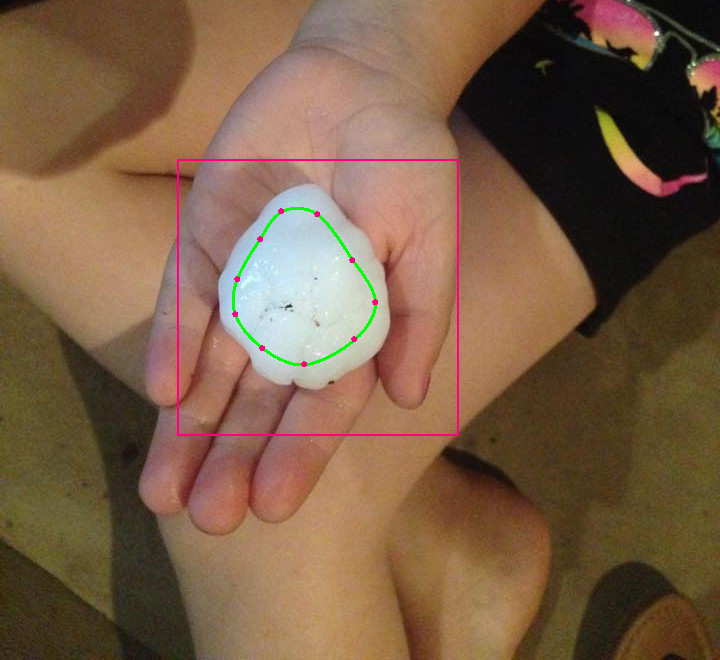}
}
\subfloat[Segmentation]{
\includegraphics[width=.24\textwidth]{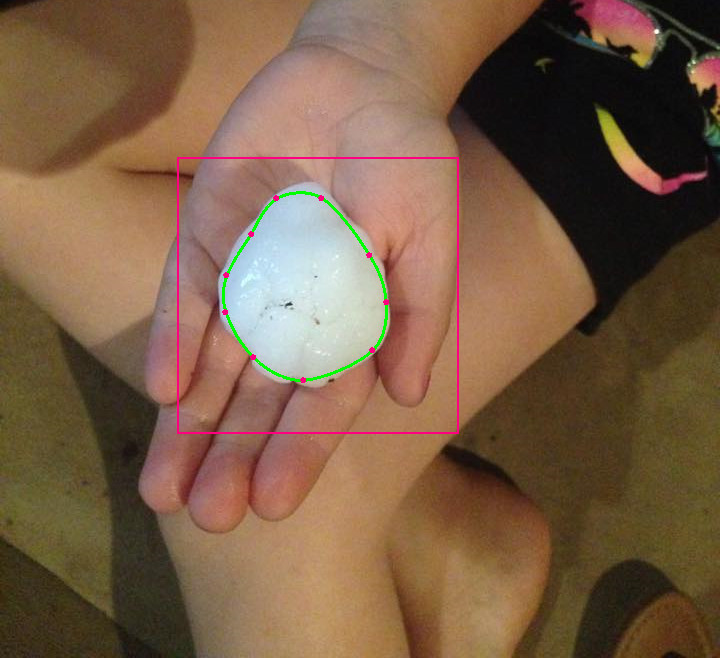}
}
\subfloat[Initialization]{
\includegraphics[width=.24\textwidth]{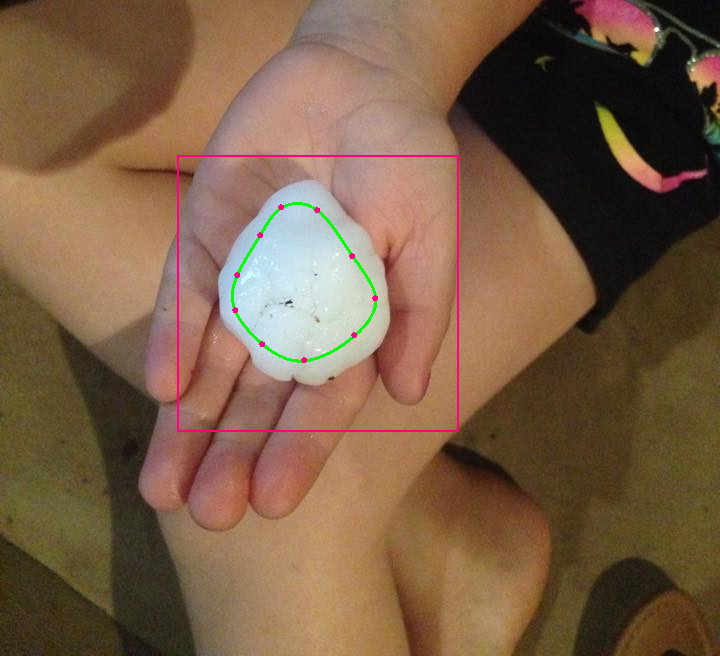}
}
\subfloat[Segmentation]{
\includegraphics[width=.24\textwidth]{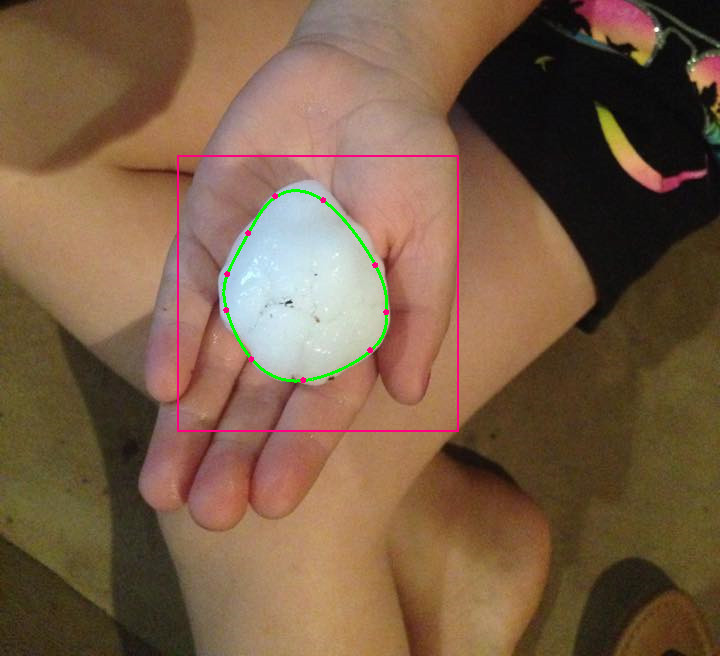}
}
\caption{Segmentation of a hail.a) and b) 4 points snake, c) and d) cubic B-spline snake.}
\label{fig:imagesHail}
\end{figure}

To show the potential of the subdivision snakes to segment objects in real images, we include in this section
the results corresponding to three images: a cerebral hemorrhage (Figure \ref{fig:imagesHemorrage}), a cell
(Figure \ref{fig:imagesCells}) and a  hail (Figure \ref{fig:imagesHail}). In general, processing real images is
more involved than  synthetic images, since several factors may affect the segmentation procedure.
In some cases the boundary of the object of interest may be blurred, see for instance Figures \ref{fig:imagesHemorrage} and \ref{fig:imagesCells}).  In other cases, the object to be segmented has inhomogeneous intensity values and poor contrast with the background, see Figure \ref{fig:imagesCells}. Moreover, sometimes  a non-uniform illumination makes difficult to capture the boundary of the object,
see Figure \ref{fig:imagesHail}. Despite these difficulties, our  method provides reasonable segmentation results, even without image enhancing in the preprocessing.

\section{Conclusions}

A method for computing the contour of an object in an image using a snake represented as a subdivision curve is presented. To illustrate its performance  we discuss the snakes associated with two classical subdivision schemes:
the four points scheme and the cubic B-spline. Our method profits from the hierarchical nature of subdivision curves, since the unknowns of the optimization process are the few control points of the subdivision curve in the coarse representation, while  good approximations of the energies and their derivatives are obtained from the fine representation.

The evolution of the snake is driven by its control points that are computed minimizing an energy which is combination of contour-based and region-based energies. We introduce a new region energy that guides the snake maximizing the contrast between the average intensity of the image within the snake and the average intensity over the complement of the snake in a fixed bounding box. Explicit expressions of the new region energy functional and its partial derivatives  are  provided and an accurate pixel discretization   is discussed.

Our experiments using synthetic and real images confirm that the proposed method is fast and robust. Our  flexible computational
framework facilitates the interaction with the snake by letting the user to move directly the control points with the mouse and
and to control the weights associated to the combination of both energy functionals. The proposed method may be extended in several
directions: the new region energy may be generalized from one channel to three channels, the subdivision schemes may be selected
among those with known formulae for the position of the limit points and their tangent vectors and different approaches of
multiresolution optimization may be applied.


\appendix

\section{Making interpolatory the cubic B-spline subdivision scheme}
\label{sec:interp_op}

The rule in \eqref{eq:limit_pos_cBS}) not only allows to evaluate the curve on dyadic parameter values, but also to impose the interpolatory condition to the scheme (see Fig. \ref{fig:BS_interpolatory}). In fact,  starting from $\bP^0$ a polygon $\widetilde{\bP}^0$  can be computed such that if,
\begin{equation}  \label{eq:operator_int}
\br(t) = \sum_{j \in \ZZ} \widetilde{\bp}_j^0 \varphi\left(t - j\right),
\end{equation}
then $\br(i)  = \bp_i^0$.

Let $\bP^0 = \left\{ \bp_i^0 \in \RR, \quad \bp_{i+M}^0 = \bp_i^0 \right\}$ be an $M$-periodic sequence of points. From \eqref{eq:operator_int} and \eqref{eq:limit_pos_cBS} it is clear that,
\begin{equation}
\bp_i^0 = \br(i) = \frac16  \widetilde{\bp}_{i-1}^0 + \frac46 \widetilde{\bp}_{i}^0 + \frac16 \widetilde{\bp}_{i+1}^0, \qquad
  \text{for} \quad  i \in \ZZ.
\end{equation}
This expression can be written in matrix terms as\footnote{The symbol $S^\infty$ comes from the limit of the subdivision operator $S$.}
\begin{equation} \label{eq:matrix_interp}
\bP^0 = S^{\infty} \widetilde{\bP}^0
\end{equation}
where $S^\infty$ is the circulant matrix with first column $b$ given by,
\begin{equation}
b = \begin{bmatrix}
{S^\infty}_{0,0} & {S^\infty}_{1,0} & \ldots & {S^\infty}_{M-1,0}
\end{bmatrix}^T
= \begin{bmatrix}
\frac23 & \frac16 & 0 & \ldots & 0 & \frac16
\end{bmatrix}^T.
\end{equation}
Hence from \eqref{eq:matrix_interp} we obtain
\[
 \widetilde{\bP} = {(S^{\infty})}^{-1} \bP^0
\]
and the problem is reduced to compute the inverse of the matrix $S^{\infty}$.

From Fourier Analysis \cite{Fraz} it is known that since $S^\infty$ is a circulant matrix, it is diagonalized by the Fourier basis. More precisely,
\begin{equation} \label{eq:FourierI_Diag}
S^\infty = \left(F^M\right)^{-1} D F^M,
\end{equation}
where $D$ is a diagonal matrix and $F^M$ is the matrix with elements,
\begin{equation} \label{eq:Fourier_basis}
F^M_{s,k} = e^{-2sk\pi\ci/M}, \qquad s,k = 0,\ldots,M-1.
\end{equation}

 It can be easily verified that,
\begin{equation*}
\left(F^M\right)^{-1} = \frac1M \overline{F^M}
\end{equation*}
where $\overline{F^M}$ is the conjugated matrix of $F^M$. Hence, from \eqref{eq:FourierI_Diag} it follows,
\begin{equation} \label{eq:Fourier_Diag}
\left(S^\infty\right)^{-1} = \frac1M \overline{F^M}\ D^{-1}\ F^M.
\end{equation}

It is also known \cite{Fraz}, that the diagonal matrix $D$ has in its diagonal the values of the Fourier transform of $b$,
\begin{equation}
\widehat{b} = \frac23 F^M_{\cdot,0} + \frac16 F^M_{\cdot,1} + \frac16 F^M_{\cdot,M-1},
\end{equation}
where $F^M_{\cdot,k}$, $k=0,\ldots,M-1$ are the column vectors with elements defined in \eqref{eq:Fourier_basis}.
Thus,
\begin{align}
\widehat{b}_s &= \frac23 F^M_{s,0} + \frac16 F^M_{s,1} + \frac16 F^M_{s,M-1}  \notag \\
	&= \frac23 + \frac16 e^{-2s\pi\ci/M} + \frac16 e^{-2(M-1)\pi\ci/M} \notag \\
	&= \frac23 + \frac13 \cos(2s\pi/M) = \widehat{b}_{M-s} \label{eq:FT_comp}
\end{align}

Since the inverse of a non-singular circulant matrix is also a circulant matrix and as $\widehat{b}_s \neq 0$ for all $s$, the matrix representing $\left(S^\infty\right)^{-1}$ is also circulant. Therefore, we only need to compute its first column $\left(S^\infty\right)^{-1}_{\cdot,0}$. From \eqref{eq:Fourier_Diag} and \eqref{eq:FT_comp} we obtain,
\begin{align}
\left(S^\infty\right)^{-1}_{s,0} &= \frac1M \sum_{t=0}^{M-1} {\widehat{b}_s}^{-1} e^{2st\pi\ci/M}  \notag \\
 &=  \begin{cases}
\frac1M + \frac3M \cos(s\pi) + \frac2M \displaystyle \sum_{t=1}^{\frac{M}{2} -1} \left( \tfrac23 + \tfrac13 \cos(2t\pi/M) \right)^{-1} \cos(2st\pi/M),  & \text{for} \quad M\mod{2} = 0 \\
\frac1M + \frac2M \displaystyle \sum_{t=1}^{\lfloor\frac{M}{2}\rfloor} \left( \tfrac23 + \tfrac13 \cos(2t\pi/M) \right)^{-1} \cos(2st\pi/M), & \text{for} \quad M\mod{2} = 1
\end{cases} \label{eq:inverse_0col}
\end{align}

\begin{remark}
To compute the entries for the whole matrix we just need to remember that it is circulant  and use \eqref{eq:inverse_0col}.
\end{remark}

\section{Computing the gradients of the energies}
\label{sec:grad_details}

With the aim to simplify the exposition, the deduction of the expressions for the gradient of the energies are showed here. Being a similar process to deduce both partial derivatives, we only show the deduction of the partial derivatives with respect to each $x_j^0$, $j=0,\ldots, M-1$.

\subsection{Gradient of {\em gradient energy}}

Deriving directly in \eqref{Egrad2} with respect to $x_j^0$ we obtain,

\begin{align}
\frac{\partial E_{grad} }{\partial x_j^0} &= - \int_{0}^M \left( \left(\frac{\partial^2 I}{\partial x^2}\frac{\partial x}{\partial x_j^0}+\frac{\partial^2 I}{\partial x \partial y}\frac{\partial y}{\partial x_j^0}\right)\frac{\dy(t)}{\dt} + \frac{\partial I}{\partial x}\frac{\partial \left(\frac{\dy}{\dt}\right) }{\partial x_j^0}\right) \dt \notag \\
 &+ \int_{0}^M \left( \left( \frac{\partial^2 I}{\partial x \partial y}\frac{\partial x}{\partial x_j^0}+ \frac{\partial^2 I}{\partial y^2}\frac{\partial y}{\partial x_j^0} \right)\frac{\dx(t)}{\dt}+ \frac{\partial I}{\partial y}\frac{\partial \left(\frac{\dx}{\dt}\right) }{\partial x_j^0} \right) \dt.
\label{parEgrad2}
\end{align}

Taking into account that $y(t)$ and $\frac{\dy(t)}{\dt}$ don't depend on $x_j^0$, from \eqref{parEgrad2} we get,
\begin{equation}
\frac{\partial E_{grad} }{\partial x_j^0} = \int_{0}^M \left(\left[\frac{\partial^2 I}{\partial x \partial y}\frac{\partial x}{\partial t}-\frac{\partial^2 I}{\partial x^2}\frac{\dy}{\dt}\right]\frac{\partial x}{\partial x_j^0}+\frac{\partial I}{\partial y}\frac{\partial \left(\frac{\dx}{\dt}\right) }{\partial x_j^0} \right) \dt
\label{derEgradxjkT}
\end{equation}

From \eqref{eq:ra_basicdef} and \eqref{eq:tanra_basicdef} it follows,
\begin{align}
\frac{\partial x(t)}{\partial x_j^0} &= \varphi(t-j), \;\;\frac{\partial y(t)}{\partial y_j^0} = \varphi(t-j) \label{deriv}\\
\frac{d x(t)}{d t} &= \sum_{j=-1}^{M+1}x_j^0 \varphi^{'}(t-j),\;\;\frac{\dy(t)}{\dt} = \sum_{j=-1}^{M+1}y_j^0 \varphi^{'}(t-j)\label{Deriv}
\end{align}
where $\varphi^{'}$ denotes $\frac{d \varphi}{dt}$. Evaluating the last expressions in $t=\frac{i}{2^k},\;\;i=0,\ldots,2^kM-1$ it holds
\begin{align}
\left.\frac{\partial x(t)}{\partial x_j^0}\right\rvert_{t=\frac{i}{2^k}} &= \varphi\left(\frac{i}{2^k}-j\right) = \left.\frac{\partial y(t)}{\partial y_j^0}\right\rvert_{t=\frac{i}{2^k}} \label{evalrx} \\
\left.\frac{\partial \left(\frac{dx}{dt}\right) }{\partial x_j^0} \right\rvert_{t=\frac{i}{2^k}} &= \varphi^{'}\left(\frac{i}{2^k}-j\right) = \left.\frac{\partial \left(\frac{dy}{dt}\right) }{\partial y_j^0} \right\rvert_{t=\frac{i}{2^k}}\label{evaldrx}
\end{align}

Substituting the integral in \eqref{derEgradxjkT} by the average of the integrand evaluated in the parameter values $i/2^k,\;\;i=0,\ldots,2^kM-1$ and using \eqref{evalrx}and \eqref{evaldrx} we obtain the following approximation for the partial derivative of gradient energy with respect to $x_j^0$,
\begin{multline}
\frac{\partial E_{grad}}{\partial x_j^0} \approx \frac{1}{2^kM}\sum_{i=0}^{2^kM-1} \left(\frac{\partial^2 I}{\partial x \partial y}\left(\br\left(\frac{i}{2^k}\right)\right) t_{ix}^k - \frac{\partial^2 I}{\partial x^2}\left(\br\left(\frac{i}{2^k}\right)\right) t_{iy}^k \right) \varphi\left(\frac{i}{2^k}-j\right) + \\ \frac{1}{2^k M}\sum_{i=0}^{2^kM-1} \frac{\partial I}{\partial y}\left(\br\left(\frac{i}{2^k}\right)\right) \varphi^{'}\left(\frac{i}{2^k}-j\right).
\end{multline}

In a similar way, deriving \eqref{Egrad2} with respect to $y_j^0$ and taking into account that $x(t)$ and $\frac{d x(t)}{dt}$ don't depend on $x_j^0$, we obtain the expression for $\frac{\partial E_{grad}}{\partial y_j^0}$.

\subsection{Gradient of {\em region energy}}

Recalling that,
\begin{equation}
\frac{\partial E_{reg}}{\partial x_j^0}=-2D\left(\frac{\partial A}{\partial x_j^0} - \frac{\partial B}{\partial x_j^0} \right)
\label{parcialEregT}
\end{equation}
where
\begin{align}
A &:= \frac{I_{\Omega}}{|\Omega|} \label{A}\\
B &:= \frac{I_R - I_{\Omega}}{\lvert R \rvert - \lvert \Omega \rvert} \label{B} \\
D &:=  A-B \label{D}.
\end{align}

Deriving directly in \eqref{A} and \eqref{B} we obtain,
\begin{align}
\frac{\partial A}{\partial x_j^0} &=\frac{1}{\lvert \Omega \rvert}\frac{\partial I_{\Omega}}{\partial x_j^0}
  - \frac{I_{\Omega}}{\lvert \Omega \rvert^2}\frac{\partial \lvert \Omega \rvert}{\partial x_j^0} \label{parcialA} \\
\frac{\partial B}{\partial x_j^0} &= \frac{-1}{\lvert R \rvert-\lvert \Omega \rvert} \frac{\partial I_{\Omega}}{\partial x_j^0}
+ \frac{I_R-I_{\Omega}}{(\lvert R \rvert-|\lvert \Omega \rvert)^2} \frac{\partial \lvert \Omega \rvert}{\partial x_j^0}.
\label{parcialB}
  \end{align}

Substituting \eqref{parcialA} and \eqref{parcialB} in \eqref{parcialEregT} we get,

\begin{equation}
\frac{\partial E_{reg}}{\partial x_j^0} = -2D\left[ \left(\frac{1}{\lvert \Omega \rvert}+\frac{1}{\lvert R \rvert-\lvert \Omega \rvert}\right)\frac{\partial I_{\Omega}}{\partial x_j^0}
-\left(\frac{I_{\Omega}}{\lvert \Omega \rvert^2}+\frac{I_R-I_{\Omega}}{\left(\lvert R \rvert-\lvert \Omega \rvert\right)^2}\right) \frac{\partial \lvert \Omega \rvert}{\partial x_j^0} \right].
\label{parcialEreg1}
\end{equation}
Now we compute the partial derivatives involved in \eqref{parcialEreg1} using the Green Theorem as stated in \eqref{Green} and \eqref{I1I2}.

Deriving in the first equality of \eqref{Green} it follows,
\[
\frac{\partial I_{\Omega}}{\partial x_j^0}=-\int_{0}^M \frac{\partial I_1}{\partial x}\frac{\partial x(t)}{\partial x_j^0}y'(t) \dt.
\]
Taking into account that, according to Leibniz's rule in \eqref{I1I2} (for differentiation under the integral sign), $\frac{\partial I_1}{\partial x}=I(x(t),y(t))$ and $\frac{\partial x(t)}{\partial x_j^0} = \varphi(t-j)$, from the previous expression we obtain,
\begin{equation}
\frac{\partial I_{\Omega}}{\partial x_j^0}=-\int_{0}^M I(\br(t))\varphi(t-j)y'(t)\;dt.
\label{derIO}
\end{equation}
Since $\lvert \Omega \rvert = \int\int_{\Omega} \dx \dy$ from \eqref{derIO} it is clear that,
\begin{equation}
\frac{\partial \lvert \Omega \rvert}{\partial x_j^0} = -\int_{0}^M \varphi(t-j)y'(t) \dt.
\label{derO}
\end{equation}

Finally, substituting \eqref{derIO} and \eqref{derO} in \eqref{parcialEreg1} and grouping similar terms we obtain
\begin{equation*}
\frac{\partial E_{reg}}{\partial x_j^0} = -2D \int_{0}^M \left[G - H\,I(\br(t)) \right]\varphi(t-j)y'(t) \dt
\end{equation*}
where
\begin{equation*}
G := \frac{I_{\Omega}}{\lvert \Omega \rvert^2}+\frac{I_R-I_{\Omega}}{(\lvert R \rvert-\lvert \Omega \rvert)^2} \quad \text{and} \quad
H := \frac{1}{\lvert \Omega \rvert}+\frac{1}{\lvert R \rvert-\lvert \Omega \rvert}.
\end{equation*}

We proceed in a similar way to compute $\dfrac{\partial E_{reg}}{\partial y_j^0}$.

\end{document}